\documentclass[11pt, reqno]{amsart}

\usepackage[usenames,dvipsnames]{color}

\usepackage[OT2, T1]{fontenc}
\usepackage{url}
\usepackage{amsmath}
\usepackage{array}
\usepackage{graphicx}
\usepackage{amssymb}
\usepackage{amsthm}
\definecolor{link}{RGB}{11,0,128}
\usepackage[colorlinks=true,citecolor=NavyBlue,linkcolor=OliveGreen,urlcolor=link]{hyperref}
\usepackage{hyperref}
\usepackage{paralist}
\usepackage{colonequals}			
\usepackage{color}
\usepackage{mathabx}			
\usepackage{amscd}
\usepackage[all,cmtip]{xy}			
\usepackage{sseq}				
\usepackage{verbatim}
\usepackage{parskip}
\usepackage{microtype}
\usepackage[in]{fullpage}
\usepackage{mathrsfs} 			
\usepackage{cleveref}
\usepackage{enumitem}
\usepackage{moreenum}			
\usepackage[alphabetic,lite]{amsrefs} 	
\usepackage{stmaryrd}			
\usepackage[foot]{amsaddr}
\usepackage{subfiles}
\usepackage{ragged2e}

\newcommand{\gA}{\alpha}

\newcommand{\bA}{\mathbb{A}}

\newcommand{\bF}{\mathbb{F}}
\newcommand{\bG}{\mathbb{G}}

\newcommand{\bL}{\mathbb{L}}

\newcommand{\bP}{\mathbb{P}}
\newcommand{\bQ}{\mathbb{Q}}

\newcommand{\bZ}{\mathbb{Z}}

\newcommand{\bbB}{\mathbf{B}}

\newcommand{\bbG}{\mathbf{G}}

\newcommand{\cA}{\mathcal{A}}
\newcommand{\cB}{\mathcal{B}}

\newcommand{\cE}{\mathcal{E}}

\newcommand{\cG}{\mathcal{G}}
\newcommand{\cH}{\mathcal{H}}

\newcommand{\cM}{\mathcal{M}}

\newcommand{\cO}{\mathcal{O}}

\newcommand{\cT}{\mathcal{T}}

\newcommand{\cW}{\mathcal{W}}
\newcommand{\cX}{\mathcal{X}}
\newcommand{\cY}{\mathcal{Y}}

\newcommand{\fm}{\mathfrak{m}}

\newcommand{\fp}{\mathfrak{p}}

\newcommand{\sB}{\mathscr{B}}

\newcommand{\sE}{\mathscr{E}}
\newcommand{\sF}{\mathscr{F}}
\newcommand{\sG}{\mathscr{G}}
\newcommand{\sH}{\mathscr{H}}
\newcommand{\sI}{\mathscr{I}}

\newcommand{\sL}{\mathscr{L}}

\newcommand{\sO}{\mathscr{O}}

\newcommand{\sR}{\mathscr{R}}

\newcommand{\sV}{\mathscr{V}}

\DeclareMathOperator{\Aut}{Aut}		
\DeclareMathOperator{\Bl}{Bl}			
\DeclareMathOperator{\Char}{char}		
\DeclareMathOperator{\Dyn}{Dyn}	
\DeclareMathOperator{\Frac}{Frac}		
\DeclareMathOperator{\GL}{GL}		
\DeclareMathOperator{\Ker}{Ker}		
\DeclareMathOperator{\Out}{Out}		
\DeclareMathOperator{\Proj}{Proj}		
\DeclareMathOperator{\rad}{rad}		
\DeclareMathOperator{\Res}{Res}		
\DeclareSymbolFont{cyrletters}{OT2}{wncyr}{m}{n}
\DeclareMathSymbol{\Sha}{\mathalpha}{cyrletters}{"58}	
\DeclareMathOperator{\Spec}{Spec}		

\newcommand{\ad}{\mathrm{ad}}			

\newcommand{\ce}{\colonequals}

\newcommand{\der}{\mathrm{der}}		

\newcommand{\eps}{\varepsilon}

\newcommand{\hra}{\hookrightarrow}
\renewcommand{\i}{^{-1}}
\newcommand{\id}{\mathrm{id}}			

\newcommand{\isomto}{\overset{\sim}{\longrightarrow}}

\newcommand{\llb}{\llbracket}			
\newcommand{\llp}{(\!(}			
	
\newcommand{\ov}{\overline}
\providecommand{\p}[1]{\left(#1\right)}

\newcommand{\ra}{\rightarrow}

\newcommand{\red}{{\mathrm{red}}}		

\newcommand{\rrb}{\rrbracket}			
\newcommand{\rrp}{)\!)}			

\newcommand{\sm}{\mathrm{sm}}			

\newcommand{\SO}{{\mathrm{SO}}}		

\newcommand{\surjects}{\twoheadrightarrow}
\newcommand{\tensor}{\otimes} 			

\newcommand{\wh}{\widehat}
\newcommand{\wt}{\widetilde}

\providecommand{\up}[1]{{\upshape(}#1{\upshape)}}
\providecommand{\uref}[1]{{\upshape\ref{#1}}}
\providecommand{\uS}{{\upshape\S}}
\providecommand{\f}[2]{\frac{#1}{#2}}

\newcommand{\bconjt}{\begin{conj-tweak}}
\newcommand{\econjt}{\end{conj-tweak}}

\renewcommand{\b}{\textbf}
\providecommand{\ucolon}{{\upshape:} }
\providecommand{\uscolon}{{\upshape;} }

\newcommand{\brems}{\begin{rems} \hfill \begin{enumerate}[label=\b{\thenumberingbase.},ref=\thenumberingbase]}

\newcommand{\erems}{\end{enumerate} \end{rems}}
\newcommand{\begs}{\begin{egs} \hfill \begin{enumerate}[label=\b{\thenumberingbase.},ref=\thenumberingbase]}

\newcommand{\eegs}{\end{enumerate} \end{egs}}
\newcommand{\m}{\item}
\newcommand{\bsm}{\begin{smallmatrix}}
\newcommand{\esm}{\end{smallmatrix}}
\newcommand{\blem}{\begin{lemma}}
\newcommand{\elem}{\end{lemma}}
\newcommand{\blemt}{\begin{lemma-tweak}}
\newcommand{\elemt}{\end{lemma-tweak}}
\newcommand{\bconj}{\begin{conj}}
\newcommand{\econj}{\end{conj}}
\newcommand{\bprob}{\begin{Problem}}
\newcommand{\eprob}{\end{Problem}}
\newcommand{\bpropt}{\begin{prop-tweak}}
\newcommand{\epropt}{\end{prop-tweak}}
\newcommand{\bthmt}{\begin{thm-tweak}}
\newcommand{\ethmt}{\end{thm-tweak}}
\newcommand{\bq}{\begin{Q}}
\newcommand{\eq}{\end{Q}}
\newcommand{\benum}{\begin{enumerate}[label={{\upshape(\alph*)}}]}
\newcommand{\benuma}{\begin{enumerate}[label={{\upshape(\arabic*)}}]}
\newcommand{\benumr}{\begin{enumerate}[label={{\upshape(\roman*)}}]}
\newcommand{\eenum}{\end{enumerate}}
\newcommand{\bitem}{\begin{itemize}}
\newcommand{\eitem}{\end{itemize}}
\newcommand{\bc}{}
\newcommand{\bd}{\begin{defn}}
\newcommand{\ed}{\end{defn}}
\newcommand{\bdt}{\begin{defn-tweak}}
\newcommand{\edt}{\end{defn-tweak}}
\newcommand{\bremt}{\begin{rem-tweak}}
\newcommand{\eremt}{\end{rem-tweak}}
\newcommand{\beg}{\begin{eg}}
\newcommand{\eeg}{\end{eg}}
\newcommand{\begt}{\begin{eg-tweak}}
\newcommand{\eegt}{\end{eg-tweak}}
\newcommand{\bcort}{\begin{cor-tweak}}
\newcommand{\ecort}{\end{cor-tweak}}
\newcommand{\bcl}{\begin{claim}}
\newcommand{\ecl}{\end{claim}}

\newcommand{\x}{\text}

\newcommand{\q}{\quad}
\providecommand{\qxq}[1]{\quad\text{#1}\quad}
\providecommand{\qx}[1]{\quad\text{#1}}

\newcommand{\qq}{\quad\quad}

\newcommand{\tst}{\textstyle}
\newcommand{\ba}{\begin{aligned}}
\newcommand{\ea}{\end{aligned}}
\newcommand{\be}{\begin{equation}}
\newcommand{\ee}{\end{equation}}
\newcommand{\bpf}{\begin{proof}}
\newcommand{\epf}{\end{proof}}
\newcommand{\bthm}{\begin{thm}}
\newcommand{\ethm}{\end{thm}}
\newcommand{\bprop}{\begin{prop}}
\newcommand{\eprop}{\end{prop}}
\newcommand{\bcor}{\begin{cor}}
\newcommand{\ecor}{\end{cor}}
\newcommand{\brem}{\begin{rem}}
\newcommand{\erem}{\end{rem}}
\newcommand*{\QED}{\hfill\ensuremath{\qed}}

\usepackage{aliascnt}
\newaliascnt{numberingbase}{subsection}
\numberwithin{equation}{numberingbase}

\newtheoremstyle{thms}{0.5em}{0pt}{\itshape}{}{\bfseries}{.}{ }{}
\theoremstyle{thms}
\newtheorem{conj}[numberingbase]{Conjecture}

\newtheorem{cor}[numberingbase]{Corollary}

\newtheorem{lemma}[numberingbase]{Lemma}

\newtheorem{prop}[numberingbase]{Proposition}

\newtheorem{Q}[numberingbase]{Question}
\newtheorem{thm}[numberingbase]{Theorem}

\newtheorem{variant}[numberingbase]{Variant}

\newtheoremstyle{claims}{0.5em}{0pt}{}{}{\itshape}{.}{ }{}
\theoremstyle{claims}
\newtheorem{claim}[equation]{Claim}

\newtheoremstyle{defs}{0.5em}{0pt}{}{}{\bfseries}{.}{ }{}
\theoremstyle{defs}
\newtheorem{defn}[numberingbase]{Definition}

\newtheorem{eg}[numberingbase]{Example}
\newtheorem*{egs}{Examples}
\newtheorem{rem}[numberingbase]{Remark}

\newtheorem*{rems}{Remarks}

\Crefname{claim}{Claim}{Claims}
\Crefname{bclaim}{Claim}{Claims}
\Crefname{sublemma}{Lemma}{Lemmas}
\Crefname{conj}{Conjecture}{Conjectures}
\Crefname{cor}{Corollary}{Corollaries}
\Crefname{defn}{Definition}{Definitions}
\Crefname{eg}{Example}{Examples}
\Crefname{prop}{Proposition}{Propositions} 
\Crefname{Q}{Question}{Questions}
\Crefname{rem}{Remark}{Remarks}
\Crefname{thm}{Theorem}{Theorems}
\Crefname{Theorem}{Theorem}{Theorems}
\Crefname{variant}{Variant}{Variants}

\theoremstyle{thms}
\newtheorem{thm-tweak}[subsection]{Theorem}
\Crefname{thm-tweak}{Theorem}{Theorems}
\newtheorem{lemma-tweak}[subsection]{Lemma}
\Crefname{lemma-tweak}{Lemma}{Lemmas}
\newtheorem{cor-tweak}[subsection]{Corollary}
\Crefname{cor-tweak}{Corollary}{Corollaries}
\newtheorem{prop-tweak}[subsection]{Proposition}
\Crefname{prop-tweak}{Proposition}{Propositions} 
\newtheorem{conj-tweak}[subsection]{Conjecture}
\Crefname{conj-tweak}{Conjecture}{Conjectures} 
\newtheorem{q-tweak}[subsection]{Question}
\Crefname{q-tweak}{Question}{Questions} 

\theoremstyle{defs}
\newtheorem{defn-tweak}[subsection]{Definition}
\Crefname{defn-tweak}{Definition}{Definitions}
\newtheorem{eg-tweak}[subsection]{Example}
\Crefname{eg-tweak}{Example}{Examples}
\newtheorem*{rems-tweak}{Remarks}
\newtheorem{rem-tweak}[subsection]{Remark}
\Crefname{rem-tweak}{Remark}{Remarks}

\newtheoremstyle{subsection-tweak}
   {2pt}
   {3pt}%
   {}
   {}%
   {\bfseries}
   {}%
   {.5em}
   {\thmnumber{\@{#1}{}\@{#2}.}%
    \thmnote{~{\bfseries#3.}}}    
    
\theoremstyle{subsection-tweak}
\newtheorem{pp}[numberingbase]{}
\newcommand{\bpp}{\begin{pp}}
\newcommand{\epp}{\end{pp}}
\newcommand{\bppt}{\begin{pp-tweak}}
\newcommand{\eppt}{\end{pp-tweak}}

\theoremstyle{subsection-tweak}
\newtheorem{pp-tweak}[subsection]{}

\makeatletter
\def\@tocline#1#2#3#4#5#6#7{
    \begingroup 
    \@ifempty{#4}{}{}

    \parindent\z@ \leftskip#3\relax \advance\leftskip\@tempdima\relax
    #5\hskip-\@tempdima
      \ifcase #1
       \or\or \hskip 2em \or \hskip 1em \else \hskip 3em \fi%
      #6\nobreak\relax
    \dotfill\hbox to\@pnumwidth{\@tocpagenum{#7}}\par
    \nobreak
    \endgroup
 }
 \def\l@section{\@tocline{1}{0pt}{1pc}{}{}}

\renewcommand{\tocsection}[3]{%
  \indentlabel{\@ifnotempty{#2}{\makebox[1.3em][l]{%
    \ignorespaces#1 \bfseries{#2}.\hfill}}}\bfseries{#3}
    \vspace{-5pt}}

\renewcommand{\tocsubsection}[3]{%
  \indentlabel{\@ifnotempty{#2}{\hspace*{-0.5em}\makebox[2.1em][l]{%
    \ignorespaces#1#2.\hfill}}}#3
    \vspace{-5pt}}
   
\makeatother 

\makeatletter 
\newcommand\appendix@section[1]{%
  \refstepcounter{section}%
  \orig@section*{Appendix \@Alph\c@section. #1}%
}
\let\orig@section\section
\g@addto@macro\appendix{\let\section\appendix@section}
\makeatother

\makeatletter
\@namedef{subjclassname@2020}{%
  \textup{2020} Mathematics Subject Classification}
\makeatother

\author{K\k{e}stutis \v{C}esnavi\v{c}ius}
\address{CNRS, Universit\'{e} Paris-Saclay,   Laboratoire de math\'{e}matiques d'Orsay, F-91405, Orsay, France}
\email{kestutis@math.u-psud.fr}

\date{\today}

\usepackage{stackengine}


\begin{document}

\title{Grothendieck--Serre in the quasi-split unramified case}

\maketitle


\begin{abstract}
The Grothendieck--Serre conjecture predicts that every generically trivial torsor under a reductive group scheme $G$ over a regular local ring $R$ is trivial. We settle it in the case when $G$ is quasi-split and $R$ is unramified. Some of the techniques that allow us to overcome obstacles that have so far kept the mixed characteristic case out of reach include a version of Noether normalization over discrete valuation rings, as well as a suitable presentation lemma for smooth relative curves in mixed characteristic that facilitates passage to the relative affine line via excision and patching.
\end{abstract}

\vspace{-50pt}

\hypersetup{
    linktoc=page,     
}

\renewcommand*\contentsname{}
\q\\
\tableofcontents

\vspace{-20pt}

\section{The Grothendieck--Serre conjecture}

The subject of this article is the following conjecture, due to Serre \cite{Ser58b}*{page~31, Remarque} and Grothendieck \cite{Gro58}*{pages~26--27, Remarques~3$^{\circ}$}, \cite{Gro68b}*{Remarques~1.11~a)}, on triviality of torsors under reductive groups.

\bconjt[Grothendieck--Serre] \label{conj:GS}
For a regular local ring $R$ and a reductive $R$-group scheme $G$, no nontrivial $G$-torsor trivializes over the fraction field of $R$, in other words,
\[
\Ker\p{H^1(R, G) \ra H^1(\Frac(R), G)} = \{*\}.
\]
\econjt

The conjecture is settled when $R$ contains a field but its remaining mixed characteristic case has so far been widely open: see the recent survey \cite{Pan18} for a detailed review of the state of the art, as well as \S\ref{pp:known-cases} below for a summary. The goal of this article is to resolve the mixed characteristic case 
under the assumption that $R$ is unramified and $G$ is quasi-split. We recall that a regular local ring $(R, \fm)$ with $p \ce \Char(R/\fm)$ is \emph{unramified} if the ring $R/pR$ is also regular, in other words, if either $R$ contains a field or $p \not\in\fm^2$. 

\begin{thm-tweak} \label{thm:main-ann}
For an unramified regular local ring $R$ and a quasi-split reductive $R$-group scheme $G$, 
\[
\Ker\p{H^1(R, G) \ra H^1(\Frac(R), G)} = \{*\};
\]
moreover, a reductive $R$-group scheme $H$ is split if and only if its generic fiber $H_{\Frac(R)}$ is split. 
\end{thm-tweak}

In fact, our result is stronger: we allow the regular ring $R$ to be semilocal, flat, and geometrically regular over some Dedekind ring $\cO$, so that the case $\cO = \bZ$ with $R$ local recovers the above, see \Cref{thm:main,thm:split} and \Cref{eg:geom-regular}. 
 The semilocal version is worth the extra effort because in many ways it is a more natural starting point. In equal characteristic, we strengthen the last aspect of \Cref{thm:main-ann}: for an equicharacteristic  regular semilocal $R$, we show that a reductive $R$-group scheme $H$ is quasi-split if and only if its generic fiber $H_{\Frac(R)}$ is quasi-split, see \Cref{thm:quasi-split}.
 





The Grothendieck--Serre conjecture is known for its numerous concrete consequences. We illustrate them with the following product formula that seems to resist any direct attack. For a further consequence that concerns quadratic forms over regular semilocal rings, see \Cref{cor:q-forms}.

\begin{cor-tweak} \label{cor:prod-formula}
For an unramified regular local ring $R$, an $r \in R \setminus \{ 0 \}$, and the $r$-adic completion $\wh{R}$,
\[
\tst G(\wh{R}[\f1r]) = G(\wh{R})G(R[\f1r]) \qxq{for every quasi-split reductive $R$-group $G$.}
\]
\end{cor-tweak}

Indeed, if the double coset on the right side did not exhaust the left side, then one could use patching (for instance, \Cref{lem:LMB-patching}) to build a nontrivial $G$-torsor that trivialized over $R[\f1r]$ (and also over~$\wh{R}$).

\begin{pp-tweak}[Known cases] \label{pp:known-cases}
Previous results on \Cref{conj:GS} fall into the following categories. 
\benuma
\m \label{KC-1}
The case when $G$ is a torus was settled by Colliot-Th\'{e}l\`{e}ne and Sansuc in \cite{CTS87}.

\m \label{KC-2}
The case when $R$ is $1$-dimensional, that is, a discrete valuation ring, was settled by Nisnevich in \cite{Nis82}, \cite{Nis84}, with corrections and a generalization to semilocal Dedekind rings by Guo in \cite{Guo20}. Subcases of the $1$-dimensional case (resp.,~of its semilocal generalization) appeared in \cite{Har67}, \cite{BB70}, \cite{BT87} (resp.,~\cite{PS16}, \cite{BvG14}, \cite{BFF17}, \cite{BFFH19}).

\m \label{KC-3}
The case when $R$ is Henselian was settled in \cite{BB70} and \cite{CTS79}*{Assertion~6.6.1}. For such $R$, one may test the triviality of a $G$-torsor after base change to the residue field, so one may choose a height $1$ prime $\fp \subset R$ for which $R/\fp$ is regular, apply the Nisnevich result, and induct on $\dim R$.  

\m \label{KC-4}
The case when $R$ contains a field, that is, when $R$ is of equicharacteristic, was settled by Fedorov--Panin in \cite{FP15} when the field is infinite (with significant inputs from \cite{PSV15}, \cite{Pan20b}), and by Panin \cite{Pan20a} when the field is finite, with substantial simplifications due to Fedorov \cite{Fed21a}. Various subcases of the equicharacteristic case appeared in \cite{Oja80}, \cite{CTO92}, \cite{Rag94},  \cite{PS97}, \cite{Zai00}, \cite{OP01}, \cite{OPZ04},  \cite{Pan05}, \cite{Zai05}, \cite{Che10}, \cite{PSV15}. 

\m \label{KC-5}
Sporadic cases, in which either $R$ or $G$ is of specific form but with $R$ possibly of mixed characteristic, appeared in \cite{Gro68b}*{Remarques~1.11~a)}, \cite{Oja82}, \cite{Nis89}, \cite{Fed21b}, \cite{Fir22}, 
\cite{BFFP20}, \cite{Pan19c}.
\eenum
In the cases \ref{KC-1}--\ref{KC-4}, the results also include the variant when the regular ring $R$ is only semilocal. 

For arguing Theorems \ref{thm:main-ann} and \ref{thm:main}, we use the toral case \cite{CTS87} and the semilocal Dedekind case \cite{Guo20} but no other known case of \Cref{conj:GS}. In fact, our argument simultaneously streamlines the case when $R$ contains a field, 
although we do not pursue this here beyond the case of quasi-split $G$ contained in \Cref{thm:main-ann} because we see no point in repeating the same additional reductions that Fedorov--Panin, Panin, and Fedorov used for handling general $G$ over such $R$.
\end{pp-tweak}

\begin{pp-tweak}[The point of departure] \label{pp:departure}
A key feature of the Grothendieck--Serre conjecture and, in fact, of problems of its flavor, for example, of the Bass--Quillen conjecture, is that one cannot easily ``enlarge'' the ring $R$, essentially, because this may trivialize torsors, one can only ``shrink'' it. The key to progress therefore lies in better understanding the geometry of $R$, and our point of departure is precisely in this for unramified $R$ of mixed characteristic $(0, p)$: we apply Popescu approximation to assume that $R$ is essentially smooth over $\bZ_{(p)}$ and then use 
a presentation theorem in the style of Gabber--Gros--Suwa \cite{CTHK97}*{Theorem~3.1.1} 
to spread $R$ out to a fibration $U \ra S$ into smooth affine curves over an open 
\[
S \subset \bA^{\dim(R) - 1}_{\bZ_{(p)}}
\]
in such a way that a given small closed subscheme $Y \subset \Spec(R)$ spreads out to be \emph{finite} over $S$, see \Cref{cor:main-structural} for a precise statement. This structural result may be viewed as a version of the Noether normalization in mixed characteristic and is reminiscent of presentation lemmas of Quillen and Gabber from \cite{Qui73}*{Lemma~5.12} and \cite{Gab94b}*{Lemma~3.1}. 

For us, $Y$ is such that a generically trivial $G$-torsor $E$ that we wish to trivialize reduces to a $B$-torsor over $\Spec(R) \setminus Y$ for a Borel $B \subset G$. The valuative criterion of properness applied to $E/B$ allows us to make this $Y$ be of codimension $\ge 2$ in $\Spec(R)$, and this codimension requirement appears difficult to relax while arguing our mixed characteristic ``Noether normalization.'' In equal characteristic, $Y$ being of codimension $\ge 1$ suffices and is immediate to arrange from the generic triviality of $E$ without using any Borel, and this distinction is one of the main sources of complications in comparison to works of Panin and Fedorov. Although in mixed characteristic virtually every step seems to require either new ideas or new techniques, the works of Panin and Fedorov in equicharacteristic have provided us with invaluable guidance for what the structure of the overall argument might be.
\end{pp-tweak}

\begin{pp-tweak}[The stages of the proof of \Cref{thm:main-ann}] \label{pp:overview}
In \Cref{thm:main-ann}, the key assertion is the triviality of every generically trivial $G$-torsor $E$. For this, our argument proceeds as follows.
\benuma
\m
In \S\ref{appendix}, we follow a suggestion by one of the referees and present a quick reduction of the quasi-split case of the Grothendieck--Serre conjecture to the setting when  our quasi-split $G$ is also semisimple and simply connected. This additional assumption simplifies the steps \ref{O-3} and \ref{O-5} below, but it could also be avoided at the cost of more refined techniques in these steps. Upon insistence of a referee, however, we take full advantage of this simplification.

\m \label{O-1}
In \S\ref{sec:fibrations}, we use the Bertini theorem and some version of the  Gabber--Gros--Suwa presentation theorem in the context of semilocal Dedekind bases $\cO$ to build the aforementioned fibration $U \ra S$ starting from a projective, flat compactification of $R$ over $\cO$, see \Cref{var:Dedekind}. We do not separate into cases according to whether the residue fields of $\cO$ are all infinite or not, but finite residue fields lead to complications that concern Bertini theorems (as do imperfect residue fields: in the case of perfect residue fields, we could replace our version of the Gabber--Gros--Suwa presentation theorem with ideas from Artin's technique of ``good neighborhoods'' from \cite{SGA4III}*{Expos\'{e}~XI}\footnote{See the previous arXiv versions of this article for detailed arguments along such lines.}). We resolve these complications 
with Gabber's approach \cite{Gab01} to Bertini theorems in positive characteristic. \emph{Op.~cit.}~is more convenient for us than the generally finer approach of Poonen \cite{Poo04} because it can guarantee that a suitable hypersurface exists for every large enough degree divisible by the characteristic, which helps in making sure that this degree is uniform across all the residue fields of $\cO$ at maximal ideals.

In this step, a major simplification in comparison to the strategy in equicharacteristic is that we do not seek a fibration into \emph{projective} curves (nor even the complement of a relatively finite subscheme in a projective relative curve) but are nevertheless able to ensure that $Y$ spreads out to a \emph{finite} $S$-scheme in the notation of \S\ref{pp:departure}. Even in equicharacteristic, this allows us to dispose of much effort usually spent in analyzing the ``boundary'' in subsequent steps.

\m \label{O-2}
In \S\ref{sec:Macaulayfication}, we deduce the mixed characteristic ``Noether normalization'' mentioned in \S\ref{pp:departure} and then use it to lift our generically trivial $G$-torsor $E$ to a torsor $\sE$ over a smooth affine $R$-curve $C$ equipped with a section $s \in C(R)$ such that $\sE$ pulls back to $E$ via $s$ and reduces to a torsor under the unipotent radical of a Borel over $C \setminus Z$ for some $R$-finite $Z \subset C$. The $R$-finiteness (as opposed to mere $R$-quasi-finiteness) of $Z$ is critical for later steps and comes from the finiteness of the spreading out of $Y$. The appearance of the unipotent radical is a new phenomenon: in equicharacteristic, $C \setminus Z$ is affine and $\sE|_{C \setminus Z}$ is a trivial torsor.

Our $(C, s, Z)$ is a simplification of what Panin and Fedorov  keep track of with the notion of a ``nice triple.'' 
 The latter is a variant of a ``standard triple'' of Voevodsky \cite{MVW06}*{Definition~11.5}  used in his construction of the triangulated category of motives. In general, it is convenient to work in terms of the relative $R$-curve $C$ instead of directly with the fibration $U \ra S$ because this gives the flexibility of changing $C$. In this, we reap the benefits of our $C$ being affine: we need to work less in subsequent reductions than ``nice triples'' would require.

\m \label{O-3}
Our $\sE$ is not a $G_C$-torsor but a $\sG$-torsor for some reductive $C$-group scheme $\sG$ equipped with a Borel $\sB \subset \sG$ whose $s$-pullback is $B \subset G$, so, in order to proceed, in \S\ref{sec:fix-G} we modify $C$ to reduce to $\sB \subset \sG$ being $B_C \subset G_C$. For this, we use the locally constant nature of the scheme parametrizing isomorphisms between two quasi-pinned reductive group schemes to show, at the cost of shrinking $C$, that $\sG$ and $G_C$ become isomorphic (compatibly with the Borels) over some \emph{finite} \'{e}tale $\wt{C} \ra C$ to which $s$ lifts.




\m \label{O-4}
After simplifying $\sG$ in \S\ref{sec:fix-G}, we turn to simplifying $C$ in \S\S\ref{sec:Lindel}--\ref{sec:to-A1}, namely, to replacing $C$ by $\bA^1_R$. In \S\ref{sec:Lindel}, we construct an affine open $U \subset C$ containing $Z$ and $s$ as well as a quasi-finite $R$-morphism $\pi\colon C \ra \bA^1_R$ that maps $Z$ isomorphically to a closed subscheme $Z' \subset \bA^1_R$ whose preimage in $U$ is precisely $Z$. The $R$-finiteness of $Z$ is critical for this, and the argument is simpler than its versions in the literature because $C$ is affine (as opposed to projective). It uses Panin's tricks with finite fields to first prepare $C$ and $Z$ for building $\pi$: when some residue fields of $R$ are finite, the initial $Z$ could have too many rational points to fit inside~$\bA^1_R$. 

Since $C$ is Cohen--Macaulay, our quasi-finite $\pi$ is necessarily flat, so the idea is to descend $\sE$ to a $G_{\bA^1_R}$-torsor via patching. We carry this out in \S\ref{sec:to-A1}: the main complication is the \emph{a priori} nontriviality of $\sE_{C\setminus Z}$, which we overcome by bootstrapping enough excision for $H^1(-, \sR_u(B))$ from  excision for quasi-coherent cohomology, see \Cref{lem:excision}. Relatedly, since $Z$ need not be principal, the patching is slightly more delicate than usual and uses \cite{MB96}.

\m \label{O-5}
The final step is the analysis of a $G_{\bA^1_R}$-torsor $\sE$ that is trivial away from an $R$-finite closed subscheme $Z' \subset \bA^1_R$. This is a problem of Bass--Quillen type that is significantly simpler when $G$ is semisimple, simply connected, and quasi-split (the general case could be approached using the geometry of the affine Grassmannian, similarly to the techniques of Fedorov from \cite{Fed21a}). 
The key inputs to this analysis are the study of the case when $R$ is a field carried out in \cite{Gil02} (or in the earlier \cite{RR84}) and the unramified nature of the Whitehead group \cite{Gil09}*{Fait 4.3,~Lemme~4.5} that builds on earlier work of Tits and uses the assumptions~on~$G$.
\eenum

Globally, our method may be viewed as a geometric reduction of the Grothendieck--Serre conjecture for $G$ over $R$ to its case for the torus $T\ce B/\sR_u(B)$ over $R$. It is tempting to expect that if $G$ is no longer quasi-split but has a parabolic subgroup $P \subset G$ with a Levi $M$, then one could find a way to reduce from $G$ to $M$. As it stands,  the sticking point in achieving this generalization is in the proof of \Cref{cor:C-and-G}: there we extend a $T$-torsor across a closed subscheme of codimension $\ge 2$ (across $Y$ in the notation of \S\ref{pp:departure}) and such extendibility fails beyond tori (although knowing how to resolve the Colliot-Th\'{e}l\`{e}ne--Sansuc purity question \cite{CTS79}*{Question~6.4} would help). 
\end{pp-tweak}


\begin{pp-tweak}[Conventions and notation] \label{conv}
As in \cite{SP}, our algebraic spaces need not be quasi-separated. For a scheme $S$ (resp.,~a ring $R$), we let $k_s$ (resp.,~$k_\fp$) denote the residue field of a point $s \in S$ (resp.,~a prime ideal $\fp \subset R$). Intersections $Y \cap Y'$ of closed subschemes $Y, Y' \subset S$ are always scheme-theoretic, and we recall from \cite{EGAIV1}*{Chapitre~0, D\'{e}finition 14.1.2} that $\dim(\emptyset) = -\infty$. We denote the (always open) $S$-smooth locus of an $S$-scheme $X$ by $X^\sm$. A scheme is \emph{Cohen--Macaulay} if it is locally Noetherian and its local rings are Cohen--Macaulay. We use the definition \cite{EGAIV1}*{Chapitre~0, D\'{e}finition 15.1.7, Section 15.2.2} of a regular sequence (so there is no condition on quotients being nonzero). A ring $\cO$ is \emph{Dedekind} if it is Noetherian, normal, and of dimension $\le 1$; by \cite{SP}*{Lemma~\href{https://stacks.math.columbia.edu/tag/034X}{034X}}, any such $\cO$ is a finite product of Dedekind domains.

We always consider \emph{right} torsors, for instance, we want sections of $G/H$ to give rise to $H$-torsors. As already seen in \S\ref{pp:overview}~\ref{O-3}--\ref{O-5}, we use scheme-theoretic notation when talking about torsors, that is, we base change the group in order to be unambiguous about what the base is; in the rare exceptions when this would make notation too cumbersome, we make sure that  no  confusion is possible. For a reductive group scheme $G$, we let $Z(G)$, $\rad(G)$, $G^\der$, and $G^\ad$ denote its center, maximal central torus, derived subgroup, and adjoint quotient (see \cite{SGA3IIInew}*{Expos\'{e}~XXII, Corollaire 4.1.7, D\'{e}finition 4.3.6, Th\'{e}or\`{e}me 6.2.1 (iv)}); for a semisimple group scheme $G$, we let $G^{\mathrm{sc}}$ denote its simply-connected cover (see \cite{Con14}*{Exercise~6.5.2~(i)} or \cite{torsors-regular}*{Proposition A.3.4}). For a parabolic group scheme $P$, we let $\sR_u(P)$ denote its unipotent radical constructed in \cite{SGA3IIInew}*{Expos\'{e}~XXVI, Proposition 1.6~(i)} (already in \cite{SGA3IIInew}*{Expos\'{e}~XXII, Proposition 5.6.9~(ii)} for a Borel).
\end{pp-tweak}

\subsection*{Acknowledgements} 
I thank the referees for insightful suggestions, especially, for suggesting the material of \Cref{appendix}. On several occasions during past years, I discussed the Grothendieck--Serre conjecture with Johannes Ansch\"{u}tz, Alexis Bouthier, Jean-Louis Colliot-Th\'{e}l\`{e}ne, Ning Guo, Roman Fedorov, Timo Richarz, and Peter Scholze, among others. I thank them for these conversations. I thank Jean-Louis Colliot-Th\'{e}l\`{e}ne, Uriya First, Ofer Gabber, Shang Li, Ivan Panin, Michael Rapoport, Timo Richarz, and Nguy\~{\^{e}}n Qu\'{\^{o}}c Th\'{\u{a}}ng for helpful remarks or correspondence.  I thank Vi\stackon[-9.3pt]{\^{e}$\mkern0.5mu$}{\d{}}n To\'{a}n H\d{o}c  for hospitality---a significant part of this project was completed during a visit there. This project received funding from the European Research Council under the European Union's Horizon 2020 research and innovation program (grant agreement No.~851146).

\section{The quasi-split case reduces to simply-connected groups} \label{appendix}


Following suggestions of one of the referees, we show that for quasi-split reductive groups the Grothendieck--Serre conjecture \ref{conj:GS} reduces to the case when the group is also semisimple and simply connected, see \Cref{prop:reduce-to-sc}. This reduction is very short and simple, and its key input is the following result from \cite{SGA3IIInew}.

\blemt \label{lem:H1-quasi-split}
For a scheme $S$, a semisimple $S$-group scheme $G$ that is either simply connected or adjoint, and a maximal $S$-torus with a Borel $S$-subgroup $T \subset B$, 
\[
T \cong \Res_{S'/S}(\bG_{m,\, S'}) \qx{for some finite \'{e}tale $S$-scheme $S'$};
\]
in particular, if $S$ is affine and semilocal, then $H^1(S, T) \cong H^1(S', \bG_m) \cong 0$.
\elemt

\bpf
The claim is a special case of \cite{SGA3IIInew}*{Expos\'{e}~XXIV, Proposition 3.13, Corollaire 3.14, Proposition 8.4}. The relevant $S'$ is the scheme of Dynkin diagrams of $G$ defined in \cite{SGA3IIInew}*{Expos\'{e}~XXIV, Section 3.2 and below}.
\epf

\bpropt \label{prop:reduce-to-sc}
For a Noetherian semilocal ring $R$ whose strict Henselization at any prime ideal is a unique factorization domain \up{for example, $R$ could be a regular semilocal ring}, the total ring of fractions $K \ce \Frac(R)$, and a quasi-split reductive $R$-group scheme $G$, if
\[
\Ker\p{H^1(R, (G^\der)^{\mathrm{sc}}) \ra H^1(K, (G^\der)^{\mathrm{sc}})} = \{ *\}, \qxq{then also} \Ker\p{H^1(R, G) \ra H^1(K, G)} = \{ *\}.
\]
\epropt

\bpf
Let $T \subset B \subset G$ be a maximal $R$-torus and an $R$-Borel subgroup of $G$ (see \cite{SGA3IIInew}*{Expos\'{e}~XXIV, Section 3.9}), and let $T' \subset B' \subset (G^\der)^{\mathrm{sc}}$ be the corresponding maximal $R$-torus and Borel $R$-subgroup of $(G^\der)^{\mathrm{sc}}$ (see \cite{SGA3IIInew}*{Expos\'{e}~XXVI, Proposition 1.19}). We have compatible isogenies
\[
\rad(G) \times (G^\der)^{\mathrm{sc}} \ra G \qxq{and}  \rad(G) \times T' \ra T
\] 
that have a common finite kernel $Z$ of multiplicative type, and \Cref{lem:H1-quasi-split} ensures that $H^1(R, T') = 0$. We will use the resulting commutative diagram of long exact cohomology sequences
\[
\xymatrix{
H^1(R, Z) \ar@{=}[d] \ar[r] & H^1(R, \rad(G)) \ar[d]^{\id \times \{*\}} \ar[r] & H^1(R, T) \ar[d] \ar[r] & H^2(R, Z) \ar@{=}[d] \\
H^1(R, Z) \ar[r] & H^1(R, \rad(G)) \times H^1(R, (G^\der)^{\mathrm{sc}}) \ar[r] & H^1(R, G) \ar[r] & H^2(R, Z),
}
\]
as well as its analogue over $K$. Thus, to proceed, we fix an $\gA \in \Ker(H^1(R, G) \ra H^1(K, G))$.

By \cite{CTS87}*{Theorem~4.3} (which is where we use the assumption on the strict Henselizations of $R$), the map $H^2(R, Z) \ra H^2(K, Z)$ is injective. Chasing the diagram above and its analogue over $K$, we therefore find that $\gA$ comes from some pair 
\[
(\beta, \gamma) \in H^1(R, \rad(G)) \times H^1(R, (G^\der)^{\mathrm{sc}})
\]
and that, by the nature of the second vertical arrow there, $\gamma|_K$ is trivial. 
The assumption on $(G^\der)^{\mathrm{sc}}$ then implies that $\gamma$ is trivial, so that $\gA$ lifts to a $\beta \in H^1(R, \rad(G))$ for which $\beta|_K$ lifts to $H^1(K, Z)$. The image of this $\beta$ in $H^1(R, T)$ lifts $\gA$ and is trivial over $K$. By \cite{CTS87}*{Theorem~4.1} (the Grothendieck--Serre conjecture \ref{conj:GS} for $T$), this image of $\beta$ is trivial, so $\gA$ is also trivial.
\epf




\section{Fibrations into smooth relative curves} \label{sec:fibrations}

We begin with the geometric part of our approach to the Grothendieck--Serre conjecture for quasi-split $G$: for a discrete valuation ring $\cO$ and a smooth, affine $\cO$-scheme $U$ of relative dimension $d > 0$ equipped with an $\cO$-fiberwise nowhere dense closed subscheme $Y \subset U$, we wish to construct a smooth morphism 
\[
\pi\colon U \ra S \subset \bA^{d - 1}_\cO
\]
whose fibers over the affine open $S$ are of dimension $1$ such that $Y \cap U$ is \emph{finite} over  $S$. Roughly, the idea is to cut $U$ by $d - 1$ suitably transversal hypersurfaces supplied by Bertini theorem and then let their defining equations be the images of the standard coordinates of $\bA^{d - 1}_\cO$. 
The actual argument given in \Cref{prop:imperfect}~\ref{imp-iv} and \Cref{var:Dedekind}~\ref{Ded-v} is more subtle because ensuring the $S$-finiteness of $Y \cap U$ is slightly more delicate.  To achieve this finiteness, we start from a projective compactification of $U$ and use a geometric presentation theorem of Gabber--Gros--Suwa \cite{CTHK97}*{Theorem~3.1.1} (if the residue fields of the maximal ideals of $\cO$ are perfect, ideas from Artin's construction of ``good neighborhoods'' in \cite{SGA4III}*{Expos\'{e}~XI, Proposition 3.3} suffice, too). 

Before turning to Bertini, we review the following avoidance lemma that we will use repeatedly.

\blemt \label{lem:avoidance}
For a ring $R$, a quasi-projective, finitely presented $R$-scheme $X$, a very $R$-ample line bundle $\sL$ on $X$, a finitely presented closed subscheme $Z \subset X$ not containing any positive-dimensional irreducible component of any $R$-fiber of $X$, and points $y_1, \dotsc, y_n \in X \setminus Z$, there is an $N_0 > 0$ such that for every $N \ge N_0$ there is an $h \in \Gamma(X, \sL^{\tensor N})$ whose vanishing scheme is a hypersurface $H \subset X$ containing $Z$ but not any $y_i$ or any positive-dimensional irreducible component of any $R$-fiber of $X$. 
\elemt

\bpf
The claim is a special case of \cite{GLL15}*{Theorem~5.1} (with definitions reviewed in \cite{GLL15}*{page~1207}). 
\epf

In the case when $R$ is a field, the following Bertini lemma allows us to impose a smoothness requirement on $ X^\sm \cap H$. Its most delicate case is when the base field is finite, in which it amounts to a mild extension of \cite{Gab01}*{Corollaries 1.6 and 1.7}, whose argument is actually our key technique.


\blemt \label{lem:Bertini}
Let $k$ be a field, let $X$ be a projective $k$-scheme of pure dimension, let $Y_1, \dotsc, Y_m, Z \subset X$ 
be closed subschemes with $Z = Z_1 \sqcup Z_0$ for some reduced $0$-dimensional $Z_0 \subset X^\sm$ all of whose residue fields are separable extensions of $k$, 
and fix a
\[
t \le \min(\dim(X), 
\dim(X) - \dim(Z)) \qx{
\up{recall from \uS\uref{conv} that $\dim(\emptyset) = -\infty$}.}
\]
For an ample line bundle $\sO_X(1)$ on $X$, there 
are hypersurfaces $H_1, \dotsc, H_t \subset X$ 
such that 
\benumr
\m \label{lem:Bertini-extra}
$H_1 \cap \dotsc \cap H_t$ is of pure dimension  $\dim(X) - t$ and contains $Z$\uscolon

\m \label{lem:Bertini-i}
$(X^\sm \setminus Z_{1}) \cap H_1 \cap \dotsc \cap H_t$ is $k$-smooth\uscolon 

\m \label{lem:Bertini-ii}
$\dim((Y_j \setminus Z) \cap \bigcap_{\ell \in L} H_\ell) \le \dim(Y_j \setminus Z) - \#L$ for $1 \le j \le m$ and $L \subset \{1, \dotsc, t\}$\uscolon


\eenum
moreover, we may iteratively choose $H_1, \dotsc, H_t$ so that, for each $i$, with $H_1, \dotsc, H_{i - 1}$ already fixed, $H_i$ may be chosen to have any sufficiently large degree divisible by the characteristic exponent of $k$. 
\elemt

We do not know how to ensure that the hypersurfaces $H_i$ in \Cref{lem:Bertini} all have the same degree.

\bpf
The pure dimension hypothesis means that all the irreducible components of $X$ have the same dimension, so \cite{EGAIV2}*{Proposition~5.2.1} ensures that $X$ is biequidimensional in the sense that the saturated chains of specializations of its points all have the same length. Similarly to \cite{macaulayfication}*{Section~4.1}, part \ref{lem:Bertini-extra} then ensures that each $X \cap H_1 \cap \dotsc \cap H_i$  inherits biequidimensionality, so is also of pure dimension. This reduces us to $t = 1$: by applying this case iteratively and at each step adjoining to the $Y_j$'s all their possible intersections with some of the already chosen $H_i$'s (to ensure \ref{lem:Bertini-ii}), we will obtain the general case. In the case $t = 1$, we fix closed points $y_1, \dotsc, y_{n} \in X\setminus Z$ that jointly meet every irreducible component of $X$ and of every $Y_j \setminus Z$. Both \ref{lem:Bertini-ii} and the dimension aspect of \ref{lem:Bertini-extra} will hold as soon as $H_1$ contains no $y_{j}$, so at the cost of requiring this we may forget about the $Y_j$.

For the rest of the argument, we begin with the case when $\Char k = 0$, in which we will use the ``classical'' Bertini theorem (one could also choose to skip this case because we will only use \Cref{lem:Bertini} for finite $k$). For this, we first claim that for every large $N > 0$ there are global sections $h_i$ of $\sO_X(N)$  whose common zero locus contains $Z$ and set-theoretically equals to it. Indeed, by repeatedly applying 
\cite{EGAIII1}*{Corollaire~2.2.4} to shrink the base locus, we first build such $h_{i'}'$ (resp.,~$h_{i''}''$) for some $N'$ (resp.,~$N''$) that is a power of $2$ (resp.,~of $3$), then express every large $N$ as $aN' + bN''$ with $a, b > 0$, and, finally, let $h_i$ be the collection of all the $h_{i'}'^{a}h_{i''}''^b$. 
By \cite{EGAIII1}*{Corollaire~2.2.4} and \cite{EGAIV4}*{Corollaire~17.15.9} (which uses the separable residue field assumption), granted that $N$ is sufficiently large, we may build another global section $h_0$ of $\sO_X(N)$ whose associated hypersurface contains $Z$ and is smooth at every point of $Z_0$. We adjoin this $h_0$ to the $h_i$ and 
then discard linear dependencies to assume that the $h_i$ are $k$-linearly independent. By \cite{EGAII}*{Proposition~4.2.3}, the $h_i$ determine a morphism 
\[
X\setminus Z  \ra \bP^{N'}_k
\]
such that the pullback of $\sO_{\bP^{N'}_k}(1)$ is our $\sO_{X\setminus Z}(N)$.  The hyperplanes in $\bP^{N'}_k$ and, compatibly, the nonzero $k$-linear combinations of the $h_i$ up to scaling are parametrized by the dual projective space. Due to the existence of a $k$-linear combination of the $h_i$ whose associated hypersurface does not contain a fixed $y_j$, a generic such hypersurface contains no $y_j$. Likewise, due to the openness of the smooth locus, the existence of a $k$-linear combination of the $h_i$ whose associated hypersurface is smooth at all the points in $Z_0$, and \cite{EGAIV3}*{Th\'{e}or\`{e}me~11.3.8~b$'$)$\Leftrightarrow$c)}, a generic such hypersurface is smooth at all the points in $Z_0$. Finally, by the Bertini theorem \cite{Jou83}*{Corollaire~6.11~2)}, the hypersurface $H$ associated to a generic $k$-linear combination of the $h_i$ is such that $(X^\sm \setminus Z) \cap H$ is $k$-smooth. 
In conclusion, since $k$ is infinite, we may choose our desired $H_1$ to be a generic such $H$. 

The remaining case when $\Char k = p > 0$ is a very minor sharpening of \cite{Gab01}*{Corollary~1.6} that is proved as there. Namely, we use the pure dimension hypothesis to apply \cite{Gab01}*{Theorem~1.1}\footnote{\emph{Loc.~cit.}~is stated in the case when the base field $k$ is finite but its proof continues to work whenever $k$ is any field of positive characteristic $p$.} with
\bitem
\m
$U$ there being our $X^\sm \setminus (Z \cup \{y_1, \dotsc, y_{n} \})$ and $\cE$ there being $\Omega^1_U$;

\m
$Z$ there being our $Z_1 \cup \{y_1, \dotsc, y_{n}\} \cup \bigcup_{z \in Z_0} \underline{\Spec}_{\sO_X}(\sO_X/\sI_{z}^2)$;

\m
$\Sigma$ there being our $Z_0 \cup \{y_1, \dotsc, y_{n}\}$;

\m
$m_0$ there being $0$; and

\m
$\sigma_0$ there being $0$ on our $Z_1$, a unit on each of our $y_1, \dotsc, y_{n}$, and a nonzero cotangent vector at $z \in Z_0$ on each of our $\underline{\Spec}_{\sO_X}(\sO_X/\sI_{z}^2)$;
\eitem
to find a finite set of closed points 
\[
F \subset X^\sm \setminus (Z \cup \{ y_1, \dotsc, y_{n} \})
\]
and, for every large $N$ divisible by $p$, a global section $h$ of $\sO_X(N)$ whose associated hypersurface contains $Z$, has a $k$-smooth intersection with $X^\sm \setminus (F \cup Z \cup \{y_1, \dotsc, y_{n} \})$, passes through every $z \in Z_0$ and is $k$-smooth there (for this we use \cite{EGAIV4}*{Corollaire~17.15.9} and the separable residue field assumption as in the characteristic $0$ case), and does not pass through any $y_{j}$. By \cite{EGAIII1}*{Corollaire~2.2.4}, if this $N$ is sufficiently large, then there is a global section $h'$ of $\sO_X(N/p)$ that vanishes on $Z \cup \{y_1, \dotsc, y_{n}\}$ and is such that $h + h'^p$ does not vanish at any point of $F$. We may then let $H_1$ be the hypersurface associated to $h + h'^p$. 
\epf

\bremt
In \ref{lem:Bertini-ii}, if $(Y_j \setminus Z) \cap \bigcap_{\ell \in L} H_\ell \neq \emptyset$, then the inequality is actually an equality because, unless the intersection is empty, cutting by $\#L$ hypersurfaces decreases dimension by at most $\#L$.
\eremt

For our purposes, the main drawback of \Cref{lem:Bertini} is its requirement that the residue fields of the points in $Z_0$ be separable over $k$: this is automatic if $k$ is perfect but cannot be removed if $k$ is imperfect because then \ref{lem:Bertini-extra} and \ref{lem:Bertini-i} cannot hold simultaneously (consider the case $t = \dim(X)$). To accommodate for imperfect $k$, we will relax \ref{lem:Bertini-extra} by no longer requiring that $Z \subset H_i$, see \Cref{prop:imperfect} for the precise statement, which uses the following review of weighted~blowups.

\begin{pp-tweak}[Weighted projective spaces]
 \label{pp:weighted}
For $w_0, \dotsc, w_d \in \bZ_{> 0}$, we consider the polynomial algebra $\bZ[t_0, \dotsc, t_d]$ that is $\bZ_{\ge 0}$-graded by declaring each $t_i$ to be of weight $w_i$ (and the constants $\bZ$ to be of weight $0$), and we let the resulting \emph{weighted projective space} be
\[
\bP_\bZ(w_0, \dotsc, w_d) \ce \Proj(\bZ[t_0, \dotsc, t_d]).
\]
We repeat the same construction over any scheme $S$ to build $\bP_S(w_0, \dotsc, w_d)$, although the latter is simply $\bP_\bZ(w_0, \dotsc, w_d) \times_\bZ S$ because the formation of $\Proj$ commutes with base change \cite{EGAII}*{Proposition~3.5.3}. We will only use weighted projective spaces when $w_0 = 1$, in which case the open subscheme of $\bP_S(1, w_1, \dotsc, w_d)$ given by $\{t_0 \neq 0\}$ is the affine space $\bA^{d}_S$ with coordinates 
\[
t_1/t_0^{w_1},\q \dotsc,\q t_d/t_0^{w_d}.
\] 
\end{pp-tweak}

\bppt[Weighted blowups] \label{pp:weighted-Bl}
For a scheme $X$, a line bundle $\sL$ on $X$, and sections 
\[
h_0 \in \Gamma(X, \sL^{\tensor w_0}),\q \dotsc,\q h_d \in \Gamma(X, \sL^{\tensor w_d}) \qxq{with} w_0, \dotsc, w_d > 0,
\]
we define the \emph{weighted blowup} of $X$ with respect to $h_0, \dotsc, h_d$ as
\[
\tst \Bl_X(h_0, \dotsc, h_d) \ce \underline{\Proj}_{\sO_X}(\sO_X[h_0, \dotsc, h_d]), \qxq{where} \sO_X[h_0, \dotsc, h_d] \subset \bigoplus_{n \ge 0} \sL^{\tensor n}
\]
is the quasi-coherent, graded $\sO_X$-subalgebra  generated by the $h_i$. The \emph{center} of this weighted blowup is the closed subscheme of $X$ cut out by  the $h_i$. By \cite{EGAII}*{Proposition 3.1.8~(i)}, the map 
\[
\Bl_X(h_0, \dotsc, h_d) \ra X \qx{is an isomorphism  away from the center.}
\]
By \S\ref{pp:weighted} and the functoriality of $\Proj$, the weighted blowup $\Bl_X(h_0, \dotsc, h_d)$ admits a morphism
\be \label{eqn:weighted-map}
\Bl_X(h_0, \dotsc, h_d) \ra \bP_\bZ(w_0, \dotsc, w_d) \qxq{determined by} t_i \mapsto h_i.
\ee
In the case when $w_0 = \dotsc = w_d$, our $\Bl_X(h_0, \dotsc, h_d)$ is identified with the usual blowup of $X$ along the closed subscheme cut out by the $h_i$: this is evident when also $\sL = \sO_X$, and the general case reduces to this one because  $\Proj$ is insensitive to twisting by line bundles \cite{EGAII}*{Proposition~3.1.8~(iii)}. 
\eppt

\bpropt \label{prop:imperfect}
Let $k$ be a field, let $X$ be a projective $k$-scheme of pure dimension $d$, let $\sO_X(1)$ be an ample line bundle on $X$, let $W \subset X^\sm$ be an open, let $x_1, \dotsc, x_n \in W$, and let $Y \subset X$ be a closed subscheme of codimension $> 0$. Upon replacing $\sO_X(1)$ by any large power, there exist nonzero
\[
h_0 \in \Gamma(X, \sO_X(1)), \ \  h_1 \in \Gamma(X, \sO_X(w_{1})), \ \ \dotsc,\ \ h_{d - 1} \in \Gamma(X, \sO_X(w_{d - 1})) \qxq{with} w_1, \dotsc, w_{d - 1} > 0,
\]
such that

\benumr
\m \label{imp-i}
the hypersurface $H_0 \ce V(h_0) \subset X$ does not contain any $x_i$\uscolon

\m \label{imp-ii}
the hypersurfaces $H_i \ce V(h_i) \subset X$ satisfy $Y \cap H_0 \cap \dotsc \cap H_{d - 1} = \emptyset$\uscolon

\m \label{imp-iii}
in the following commutative diagram with vertical  maps determined by the $h_0, \dotsc, h_{d - 1}$\ucolon
\[
\qq \xymatrix{
X \setminus H_0 \ar[d]_\pi \ar@{^(->}[r]  & X \setminus (H_0 \cap \dotsc \cap H_{d - 1}) \ar[d]_{\ov{\pi}} \ar@{^(->}[r] & \ov{X} \ce \Bl_X(h_0, \dotsc, h_{d - 1}) \ar[d]_{\ov{\pi}}^{\eqref{eqn:weighted-map}} \\ 
\bA^{d -1}_k \ar@{^(->}[r]^-{\x{\uS\uref{pp:weighted}}} & \bP_k(1, w_1, \dotsc, w_{d - 1}) \ar@{=}[r] & \bP_k(1, w_1, \dotsc, w_{d - 1}),
}
\]
for every $i$, the map $\pi$ is smooth \up{of relative dimension $1$} at $x_i$\uscolon

\m \label{imp-new}
for every $i$, we have $Y \cap H_0 \cap \ov{\pi}\i(\pi(x_i)) = \emptyset$\uscolon

\m \label{imp-v}
if $Y \setminus X^\sm$ is of codimension $> 1$ in $X$, then $\pi$ is smooth at each $Y \cap \ov{\pi}\i(\pi(x_i))$\uscolon

\m \label{imp-new-2}
if even $Y \setminus W$ is of codimension $> 1$ in $X$, then, for every $i$, we have $(Y \setminus W) \cap \ov{\pi}\i(\pi(x_i)) = \emptyset$\uscolon

\m \label{imp-iv}
if even $Y \setminus W$ is of codimension $> 1$ in $X$, then there are affine opens 
\[
\qq S \subset \bA^{d - 1}_k \qxq{and} x_1, \dotsc, x_n \in U \subset W \cap \pi\i(S)  \subset  X\setminus H_0
\]
such that $\pi\colon U \ra S$ is smooth of relative dimension $1$ and $Y \cap U = Y \cap \pi\i(S)$ 
is $S$-finite\uscolon
\eenum

moreover, we may iteratively choose $h_0, \dotsc, h_{d - 1}$ so that, for each $i$, with $h_0, \dotsc, h_{i - 1}$ already fixed, the degree $w_i$ may be any sufficiently large integer divisible by the characteristic exponent of $k$. 
\epropt


\bpf
Closed points are dense in a finite type scheme over a field \cite{SP}*{Lemma~\href{https://stacks.math.columbia.edu/tag/02J6}{02J6}}, so each $x_i$ specializes to a closed point of $X^\sm$. 
By replacing the $x_i$ by such specializations, we assume for the rest of the proof that the points $x_1, \dotsc, x_n$ are closed in $X$. 

In \ref{imp-iii}, the parenthetical relative dimension aspect follows from a dimension count. The intersections in \ref{imp-new} make sense because \ref{imp-ii} ensures that 
\[
Y \subset X \setminus (H_0 \cap \dotsc \cap H_{d - 1}),
\]
that is, that $Y$ does not meet the center of the weighted blowup $\ov{X} \ra X$ (see \S\ref{pp:weighted-Bl}). Moreover, \ref{imp-new} implies that over the open neighborhood of $\bigcup_{i = 1}^n \pi(x_i)$ given by the complement of $\ov{\pi}(Y \cap H_0)$, our closed $Y$ even lies in $X \setminus H_0$. Since $H_0$ is a hyperplane section, intersecting with it cuts the dimension of any positive-dimensional closed subvariety of $X$ by at most $1$, so \ref{imp-new} implies that $Y \cap \ov{\pi}\i(\pi(x_i))$ is finite. The openness of the quasi-finite locus \cite{SP}*{Lemma~\href{https://stacks.math.columbia.edu/tag/01TI}{01TI}} and the finiteness of proper, quasi-finite morphisms \cite{SP}*{Lemma~\href{https://stacks.math.columbia.edu/tag/02OG}{02OG}} then imply that $Y \cap \pi\i(S)$ is $S$-finite for every sufficiently small affine open $S \subset \bA^{d - 1}_k$ containing $\pi(x_1), \dotsc, \pi(x_n)$. Consequently, by spreading out from the semilocalization at $\pi(x_1), \dotsc, \pi(x_n)$, for any affine open $U \subset \pi\i(S)$ containing all the $x_i$ and all the points of $Y \cap \ov{\pi}(\pi\i(x_i))$, at the cost shrinking $S$ further we may also arrange that 
\[
Y \cap U = Y \cap \pi\i(S).
\]
In conclusion, \ref{imp-iv} follows from \ref{imp-i}--\ref{imp-v} (we stated \ref{imp-iv} explicitly because we will use it below).

For the rest, we may assume that $d > 0$ and begin with the case when $k$ is perfect, in which we will deduce the claim from \Cref{lem:Bertini}, even with the additional requirement that each $H_i$ passes through $x_1, \dotsc, x_n$: this will ensure that $\pi(x_i)$ is the origin of $\bA^{d - 1}_k$ for every $i$, so that 
\be \label{eqn:zero-fiber}
\ov{\pi}\i(\pi(x_i)) = (H_1 \cap \dotsc \cap H_{d - 1})\setminus (H_0 \cap \dotsc \cap H_{d - 1}) \qxq{in} X \setminus (H_0 \cap \dotsc \cap H_{d - 1}),
\ee
to the effect that \ref{imp-new} will follow from \ref{imp-ii}. We begin by using \Cref{lem:avoidance} to choose a hypersurface $H_0 \subset X$ of any sufficiently large degree $w_0$ such that $H_0$ does not contain any $x_i$ nor any irreducible component of $Y$. We then replace $\sO_X(1)$ by its $w_0$-th power and, since the extensions $k(x_i)/k$ are separable, use \Cref{lem:Bertini} to choose hypersurfaces $H_1, \dotsc, H_{d - 1}$ containing every $x_i$ such that 
\bitem
\m
$X^\sm \cap H_1 \cap \dotsc \cap H_{d - 1}$ is $k$-smooth of pure dimension $1$;

\m
$Y \cap H_0 \cap \dotsc \cap H_{d - 1} = \emptyset$;

\m
if $Y\setminus X^\sm$ is of codimension $> 1$ in $X$, then $(Y \setminus X^\sm) \cap H_1 \cap \dotsc \cap H_{d - 1} = \emptyset$; and

\m
if $Y\setminus W$ is of codimension $> 1$ in $X$, then $(Y \setminus W) \cap H_1 \cap \dotsc \cap H_{d - 1} = \emptyset$.
\eitem
These properties ensure the desired \ref{imp-i}--\ref{imp-new-2}: indeed, by \eqref{eqn:zero-fiber}, the $x_i$ and, if $Y\setminus X^\sm$ is of codimension $> 1$ in $X$, then also the $Y \cap \ov{\pi}\i(\pi(x_i))$ all lie in the smooth locus of the zero fiber of $\pi$ and, by the dimensional flatness criterion \cite{EGAIV2}*{Proposition~6.1.5}, also in the smooth locus of $\pi$. The aspect about arranging the degrees $w_i$ to be large follows from its counterpart in \Cref{lem:Bertini}.

In the remaining case when $d > 0$ and $k$ is infinite (or even imperfect), we will deduce the claim from a geometric presentation theorem of Gabber--Gros--Suwa \cite{CTHK97}*{Theorem~3.1.1}, according to which:
\be \label{eqn:CTHK-input} \tag{$\bigstar$}
\parbox{14cm}{for a smooth, affine, irreducible $k$-scheme $V$ of finite type and dimension $d$, points $v_1, \dotsc, v_m \in V$, and a closed subscheme  $Z \subset V$ of codimension $> 0$, there is a $k$-map $V \ra \bA^{d -1}_k$ that is smooth at every $v_i$ and such that $Z$ is finite over $\bA^{d - 1}_k$; without the smoothness assumption on $V$, the same holds granted that $m = 0$.}
\ee
\emph{Loc.~cit.}~is much more general but \eqref{eqn:CTHK-input} will suffice for us. For convenience, we added the aspect about nonsmooth $V$, which follows from Noether normalization \cite{SP}*{Lemma~\href{https://stacks.math.columbia.edu/tag/00OY}{00OY}}. To apply \eqref{eqn:CTHK-input}, we first use \Cref{lem:avoidance} to choose a hypersurface $H_0 = V(h_0) \subset X$ of any sufficiently large degree $w_0$ such that  $H_0$ does not contain any $x_i$, does not contain any irreducible component of $Y$ of codimension $1$ in $X$, contains the intersections of any two distinct irreducible components of $X$, contains the nonsmooth locus of any generically smooth irreducible component of $X$, and contains $Y \setminus X^\sm$ if the latter is of codimension $> 1$ in $X$. Then $X \setminus H_0$ is affine and a nonempty disjoint union of irreducible components $V$ of dimension $d$, each of which is either smooth or contains no $x_i$. By applying \eqref{eqn:CTHK-input} to each such $V$ with the $v_j$ being those $x_i$ that lie on $V$ and $Z$ being $Y \cap V$, we obtain a $k$-morphism 
\[
\pi' \colon X \setminus H_0 \ra \bA^{d - 1}_k
\]
that is smooth at every $x_i$ and such that $Y \setminus H_0$ is finite over $\bA^{d - 1}_k$. 

We replace $\sO_X(1)$ by its $w_0$-th power, and consequently view $h_0$ as a global section of $\sO_X(1)$, in other words, as a morphism $\sO_X \ra \sO_X(1)$ determined by $1 \mapsto h_0$, which for any $w > 0$ gives rise to a morphism $\sO_X \ra \sO_X(w)$ determined by $1 \mapsto h_0 \tensor \dotsb \tensor h_0$. Via the latter, any section of $\sO_X$ over $X \setminus H_0$ extends to a section of $\sO_X(w)$ over $X$ granted that $w$ is large enough, and two fixed such extensions become equal granted that we enlarge $w$ further (the uniqueness allows one to work on a finite cover of $X$ by affines that trivialize $\sO_X(1)$, and the claim then reduces to the observation that for any ring $R$, any $r \in R$, and any $h \in R[\f1r]$, the element $r^wh \in R[\f1r]$ lifts to $R$ granted that $w > 0$ is large, and any two such lifts agree granted that we multiply them by $r^{w'}$ for a large $w' > 0$). In particular, the images under $\pi'$ of the standard coordinates of $\bA^{d - 1}_k$ extend to nonzero sections 
\be \label{eqn:h-prime}
h'_1 \in \Gamma(X, \sO_X(w_1)), \q \dotsc, \q h'_{d - 1} \in \Gamma(X, \sO_X(w_{d - 1}))
\ee
for all large enough $w_1, \dotsc, w_{d - 1} > 0$. By construction, with $H_i' \ce V(h_i')$,  the map 
\[
\ov{\pi}' \colon X \setminus (H_0 \cap H_1' \cap \dotsb \cap H'_{d - 1}) \ra \bP_k(1, w_1, \dotsc, w_{d - 1})
\]
that results from $h_0$ and these $h_i'$ as in \ref{imp-iii} restricts to $\pi'$, satisfies \ref{imp-i}, and is smooth at every $x_i$. 

To proceed further, for $h_0$ and $w_1, \dotsc, w_{d - 1}$ fixed above, we consider the $k$-scheme $\cM$ that parametrizes all $(d-1)$-tuples of global sections like in \eqref{eqn:h-prime}, so that $\cM$ is noncanonically isomorphic to some affine space $\bA^N_k$. The base change $X_\cM$ comes equipped with the universal hypersurfaces 
\[
\cH_1, \dotsc, \cH_{d - 1} \subset X_\cM
\]
and, of course, also with the base changes of $x_1, \dotsc, x_n$, $Y$, and $H_0$. The defining equations of $H_0$ and of the $\cH_1, \dotsc, \cH_{d - 1}$  give the universal diagram of maps
\[
\qq \xymatrix{
X_\cM \setminus (H_0)_\cM \ar[d] \ar@{^(->}[r]  & X_\cM \setminus ((H_0)_\cM \cap \cH_1 \cap \dotsc \cap \cH_{d - 1}) \ar[d]  \\ 
\bA^{d -1}_\cM \ar@{^(->}[r]^-{\x{\uS\uref{pp:weighted}}} & \bP_\cM(1, w_1, \dotsc, w_{d - 1}) 
}
\]
whose base changes to $k$-points of $\cM$ are left squares of diagrams as in \ref{imp-iii}. To finish the proof, we will argue that for all sufficiently large $w_1, \dotsc, w_{d - 1}$ that are divisible by the characteristic exponent of $k$ and that may be chosen iteratively as in the statement, the conditions \ref{imp-i}--\ref{imp-new-2} define a nonempty open of $\cM$: since $\cM \simeq \bA^N_k$ and $k$ is infinite, this open will have a $k$-point that, by construction, will correspond to the desired $h_1, \dotsc, h_{d - 1}$ satisfying \ref{imp-i}--\ref{imp-new-2}. 

First of all, since $h_0$ is fixed, \ref{imp-i} holds by construction. Since the image of a proper morphism is closed, the nonvanishing requirement on the $h_i$ and  the condition \ref{imp-ii} define an open of $\cM$ that, granted that the $w_1, \dotsc, w_{d - 1}$ are large enough, is nonempty by \Cref{lem:avoidance}. Due to the existence of the $h_i'$ as in \eqref{eqn:h-prime} and the openness of the smooth locus, the smoothness requirement \ref{imp-iii} likewise defines a nonempty open. Granted \ref{imp-ii}, since the image of a proper morphism is closed, the requirements 
\[
Y \cap H_0 \cap \ov{\pi}\i(\pi(x_i)) = \emptyset
\]
of \ref{imp-iv} define an open. Granted that the $w_1, \dotsc, w_{d - 1}$ are iteratively chosen to be large enough and divisible by the characteristic exponent of $k$, this open is nonempty thanks to \Cref{lem:Bertini} (applied after base change to $\ov{k}$ and similarly to the case of a perfect $k$ treated above), and even remains a  nonempty open after imposing the further requirements \ref{imp-v} and \ref{imp-new-2}. In summary, the conditions \ref{imp-i}--\ref{imp-new-2} indeed define a nonempty open of $\cM$, as desired.
\epf

For our mixed characteristic purposes, we need to extend \Cref{prop:imperfect} by allowing $k$ be a semilocal Dedekind ring $\cO$. We deduce such an extension in the following variant: this  is straight-forward by specialization if the $x_i$ all lie above the maximal ideals of $\cO$, and we reduce to this case via a trick that adds an auxiliary maximal ideal to $\cO$. 

\begin{variant} \label{var:Dedekind}
Let $\cO$ be a semilocal Dedekind ring, let $X$ be a projective, flat $\cO$-scheme with fibers of pure dimension $d$, let $\sO_X(1)$ be an $\cO$-ample line bundle on $X$, let $W \subset X^\sm$ be an open, let $x_1, \dotsc, x_n \in W$, and let $Y \subset X$ be a closed subscheme that is $\cO$-fiberwise of codimension $> 0$. Upon replacing $\sO_X(1)$ by any large power, there exist nonzero
\[
h_0 \in \Gamma(X, \sO_X(1)), \ \  h_1 \in \Gamma(X, \sO_X(w_{1})), \ \ \dotsc,\ \ h_{d - 1} \in \Gamma(X, \sO_X(w_{d - 1})) \qxq{with} w_1, \dotsc, w_{d - 1} > 0,
\]
such that

\benumr
\m \label{Ded-i}
the hypersurface $H_0 \ce V(h_0) \subset X$ does not contain any $x_i$\uscolon

\m \label{Ded-ii}
the hypersurfaces $H_i \ce V(h_i) \subset X$ satisfy $Y \cap H_0 \cap \dotsc \cap H_{d - 1} = \emptyset$\uscolon

\m \label{Ded-iii}
in the following commutative diagram with vertical  maps determined by the $h_0, \dotsc, h_{d - 1}$\ucolon
\[
\qq \xymatrix{
X \setminus H_0 \ar[d]_\pi \ar@{^(->}[r]  & X \setminus (H_0 \cap \dotsc \cap H_{d - 1}) \ar[d]_{\ov{\pi}} \ar@{^(->}[r] & \ov{X} \ce \Bl_X(h_0, \dotsc, h_{d - 1}) \ar[d]_{\ov{\pi}}^{\eqref{eqn:weighted-map}} \\ 
\bA^{d -1}_\cO \ar@{^(->}[r]^-{\x{\uS\uref{pp:weighted}}} & \bP_\cO(1, w_1, \dotsc, w_{d - 1}) \ar@{=}[r] & \bP_\cO(1, w_1, \dotsc, w_{d - 1}),
}
\]
for every $i$, the map $\pi$ is smooth \up{of relative dimension $1$} at $x_i$\uscolon

\m \label{Ded-new}
for every $i$, we have $Y \cap H_0 \cap \ov{\pi}\i(\pi(x_i)) = \emptyset$\uscolon

\m \label{Ded-iv}
if $Y \setminus X^\sm$ is $\cO$-fiberwise of codimension $> 1$ in $X$, then $\pi$ is smooth at each $Y \cap \ov{\pi}\i(\pi(x_i))$\uscolon

\m \label{Ded-new-2}
if even $Y \setminus W$ is $\cO$-fiberwise of codimension $> 1$ in $X$, then $(Y \setminus W) \cap \ov{\pi}\i(\pi(x_i)) = \emptyset$ for~all~$i$\uscolon

\m \label{Ded-v}
if even $Y \setminus W$ is $\cO$-fiberwise of codimension $> 1$ in $X$, then there are affine opens 
\[
\qq S \subset \bA^{d - 1}_\cO \qxq{and} x_1, \dotsc, x_n \in U \subset W \cap \pi\i(S) \subset  X\setminus H_0
\]
such that $\pi\colon U \ra S$ is smooth of relative dimension $1$ and $Y \cap U = Y \cap \pi\i(S)$ 
is $S$-finite.
\eenum
\end{variant}

\bpf
If $\cO$ is a field, then the claim follows from \Cref{prop:imperfect}, so we may assume that $\cO$ is a domain that is not a field and set $K \ce \Frac(\cO)$. As in the proof \Cref{prop:imperfect}, we may replace each $x_i$ by its suitable specialization to assume without losing generality that the points $x_1, \dotsc, x_n$ are closed in their $\cO$-fibers of $X$. Moreover, \ref{Ded-v} again follows from the rest by an argument as~there. 

For the rest, we let $C \subset \Spec(\cO)$ be the union of the closed points of $\Spec(\cO)$, and we begin with the case when each $x_i$ lies over $C$. By \Cref{prop:imperfect}, the desired $h_i$ exist after base change to $C$. Moreover, the last aspect of \Cref{prop:imperfect} ensures that these $h_i$ may be chosen to have constant degrees $w_i$ on $C$ and, by \cite{EGAIII1}*{Corollaire~2.2.4}, to be such that they lift to sections of $\sO_X(w_i)$. Since the image of a proper morphism is closed, these lifts still satisfy the desired \ref{Ded-i}--\ref{Ded-new-2} (the formation of the weighted blowup $\Bl_X(h_0, \dotsc, h_{d - 1})$ in \ref{Ded-iii} need not commute with base change to $C$, but the formation of the map 
\[
\ov{\pi} \colon X \setminus (H_0 \cap \dotsc \cap H_{ d - 1}) \ra \bP_\cO(1, w_1, \dotsc, w_{d - 1})
\]
does and this suffices).

More generally, consider the case when each $x_i$ that lies in the generic $\cO$-fiber of $X$ specializes to some point in $W \cap X_C$. By replacing each such $x_i$ by this specialization, we reduce to $x_1, \dotsc, x_n$ all lying over $C$, that is, to the previous paragraph. In general, such specializations need not exist, but we will reduce to the case when they do by enlarging $\Spec(\cO)$ and suitably extending $X$ and $Y$. For this, we begin by letting $\fm \subset \cO$ range over the maximal ideals and noting that, since $\bigcup_\fm \Spec(\cO_\fm)$ is an open cover of $\Spec(\cO)$,  the map 
\[
\tst \cO \isomto \bigcap_\fm \cO_\fm
\]
is an isomorphism. For every subfield $K' \subset K$, we consider the discrete valuation rings $\cO'_\fm \ce \cO_\fm \cap K'$ and we note that, by \cite{Mat89}*{Theorem~12.2}, the intersection 
 \[
 \tst \cO' \ce \bigcap_\fm \cO_\fm'
 \] 
 is a semilocal Dedekind subdomain of $\cO$ whose localizations at maximal ideals are precisely the $\cO'_\fm$. As $K'$ ranges over the finitely generated subfields of $K$, the subrings $\cO'$ exhaust $\cO$. Thus, we may use a limit argument to replace $\cO$ by some such $\cO'$, and hence to reduce to the case when $K$ is finitely generated (and automatically separable) over its prime subfield. By \cite{EGAIV4}*{Corollaire~17.15.9} and spreading out, such a $K$ is the fraction field of a regular domain $A$ that is smooth either over this prime subfield or even over $\bZ$. Moreover, this $A$ cannot be $0$-dimensional: else, $K$ would be a finite field, which would contradict our assumption that $\cO$ is not a field. 

By localizing $A$ if needed, we may assume that 
\bitem
\m
$X_K$ spreads out to a projective, flat $A$-scheme $\cX$ with $A$-fibers of pure dimension $d$ (see \cite{EGAIV3}*{Th\'{e}or\`{e}me~12.2.1 (ii) and (v)});

\m
$\sO_{X_K}(1)$ spreads out to an $A$-ample line bundle $\sL$ on $\cX$ (see \cite{EGAIV3}*{Corollaire~9.6.4});

\m
$W_K$ and $Y_K$ spread out to $A$-flat open (resp.,~closed) subschemes 
\[
\qq \cW \subset \cX^\sm \qxq{and} \cY \subset \cX
\]
such that $\cY$ is $A$-fiberwise of codimension $> 0$ and, if $Y_K\setminus W_K$ (resp.,~$Y_K \setminus X^\sm_K$) is of codimension $> 1$ in $X$, then $\cY \setminus \cW$ (resp.,~$\cY \setminus \cX^\sm$) is $A$-fiberwise of codimension $> 1$ in $\cX$ (see \cite{EGAIV3}*{Th\'{e}or\`{e}me~11.3.1, Corollaire 12.2.2~(i)}); and 

\m
each $x_i$ that lies over $K$ spreads out to an $A$-finite closed subscheme of $\cW$.
\eitem
 Since $A$ is positive-dimensional, it has infinitely many primes $\fp$ of height $1$. In particular, for some such $\fp$, the discrete valuation subring $R \ce A_\fp$ of $K$ is distinct from each $\cO_\fm$. Since these discrete valuation subrings share a common fraction field, there can be no inclusion relations between them. Thus, \cite{Mat89}*{Theorems~12.2 and~12.6} ensure that 
 \[
 \tst \wt{\cO} \ce R \cap \bigcap_\fm \cO_\fm
 \]
 is a semilocal Dedekind domain whose spectrum is obtained by glueing $\Spec(\cO)$ and $\Spec(R)$ along their common open $\Spec(K)$. Consequently, we may glue $X$ with $\cX_R$ along $X_K$ to extend $X$ to a projective, flat $\wt{\cO}$-scheme $\wt{X}$ with fibers of pure dimension $d$ that comes equipped with an ample line bundle $\sO_{\wt{X}}(1)$ obtained by glueing $\sO_X(1)$ with $\sL_R$ along $\sO_{X_K}(1)$, with an open subscheme $\wt{W} \subset \wt{X}^\sm$ obtained by glueing $W$ with $\cW_R$ along $W_K$, as well as with a closed subscheme $\wt{Y} \subset \wt{X}$ obtained by glueing $Y$ with $\cY_R$ along $Y_K$. We may now replace $X$ and $\cO$ by $\wt{X}$ and $\wt{\cO}$ without losing generality. This has the advantage that, by construction, we are left with the already treated case in which each $x_i$ that lies in $X_K$ specializes to some point in $W \cap X_C$. 
\epf

\section{Lifting the torsor to a smooth relative curve}\label{sec:Macaulayfication}

In practice, we start from a smooth affine scheme $W$, not from its projective compactification $X$ as in \Cref{var:Dedekind}. The following proposition recasts the results of the previous section from this vantage point. We thank Panin for extracting its formulation from the initial version of this article.

\bpropt \label{cor:main-structural}
For 
\bitem
\m
a semilocal Dedekind ring $\cO$\uscolon

\m
a smooth $\cO$-algebra $A$ that is everywhere of positive relative dimension over $\cO$\uscolon

\m
$x_1, \dotsc, x_n \in \Spec(A)$\uscolon 

\m
a closed subscheme $Y \subset \Spec(A)$ that is of codimension $\ge 2$\uscolon
\eitem
there are
\benumr
\m \label{MS-i}
an affine open $U \subset \Spec(A)$ containing all the $x_i$\uscolon

\m
an affine open $S \subset \bigsqcup_{d \ge 0} \bA^d_\cO$\uscolon and

\m \label{MS-iii}
a smooth $\cO$-morphism $\pi\colon U \ra S$ of pure relative dimension $1$ such that $Y \cap U$ is $S$-finite.
\eenum
\epropt

\bpf
We decompose $\Spec(\cO)$ and $W \ce \Spec(A)$ into connected components to assume that $\cO$ and $A$ are domains, so that $A$ is of pure relative dimension $d > 0$ over $\cO$. 
We embed $W$ into some affine space over $\cO$ and then take the schematic image in the corresponding projective space to build an open immersion $W \hra X$ into a projective, flat $\cO$-scheme $X$. On the $K$-fiber with $K \ce \Frac(\cO)$ this immersion has dense image, so $X_K$ is of pure dimension $d$. It then follows from \cite{SP}*{Lemmas~\href{https://stacks.math.columbia.edu/tag/0D4J}{0D4J} and \href{https://stacks.math.columbia.edu/tag/02FZ}{02FZ}} that $X$ is of pure relative dimension $d$ over $\cO$. We will obtain our $U$ and $\pi$ from $X$ via \Cref{var:Dedekind}~\ref{Ded-v}. To apply the latter, all we need to do is to check that the schematic image $\ov{Y}$ of $Y$ in $X$ is such that $Y' \ce \ov{Y} \setminus Y$ is $\cO$-fiberwise of codimension $> 1$ in $X$.


By \cite{SP}*{Lemma \href{https://stacks.math.columbia.edu/tag/01R8}{01R8}}, set-theoretically we have $\ov{Y} = \bigcup_{y} \ov{\{ y\}}$ where  $y$ ranges over the generic points of $Y$ and the closures are in $X$. Each $y$ is of height $\ge 2$ in $X$, so each $\ov{\{y \}}$ intersects the $\cO$-fiber of $X$ containing $y$ in a closed subscheme of dimension $\le d - 1$ (even $\le d - 2$ if the fiber is generic). Thus, since $\ov{\{y\}}$ has a nonempty open $\ov{\{y \}} \cap W$, the contribution of $y$ to its $\cO$-fiber of $Y'$ is of dimension $\le d - 2$. The only situation in which $\ov{\{y\}}$ may contribute to other $\cO$-fibers of $Y'$ is when $y$ lies in the generic $\cO$-fiber of $X$ and $\cO$ is not a field. However, since the local rings of $X$ are of dimension $\le d + 1$, then the intersection of $\ov{\{y\}}$ with any closed $\cO$-fiber of $X$ is of dimension $\le d - 2$. In conclusion, $Y'$ is $\cO$-fiberwise of dimension $\le d - 2$, that is, $\cO$-fiberwise of codimension $> 1$ in $X$. 
\epf


The following consequence of \Cref{cor:main-structural} starts a string of reductions that will eventually lead to \Cref{thm:main-ann}. In comparison to versions in the literature, for instance, to \cite{Fed21b}*{Proposition~4.4}, the main new phenomena are that the group $\sG$ is only defined over a small affine $C$ and that $C\setminus Z$ need not be affine, which will cause additional subtleties in \S\ref{sec:to-A1}. It is not necessary to assume that our quasi-split group is simply connected but, upon insistence of a referee, we assume this anyway.  




\bpropt \label{cor:C-and-G}
For a semilocal Dedekind ring $\cO$, the semilocalization $R$ of a smooth $\cO$-algebra $A$ at finitely many primes $\fp$, 
a quasi-split reductive $R$-group $G$, a Borel $R$-subgroup $ B \subset G$, 
and a generically trivial $G$-torsor $E$,
there are
\benumr


\m \label{CG-D}
a smooth, affine $R$-scheme $C$ of pure relative dimension $1$\uscolon


\m \label{CG-iii}
a section $s \in C(R)$\uscolon

\m \label{CG-iv}
an $R$-finite closed subscheme $Z \subset C$\uscolon

\m \label{CG-v-new}
a quasi-split reductive $C$-group scheme $\sG$ with a Borel $\sB \subset \sG$ whose $s$-pullback is $B\subset G$\uscolon

 \m \label{CG-v}
a $\sG$-torsor $\sE$ whose $s$-pullback is $E$ such that $\sE$ reduces to an $\sR_u(\sB)$-torsor over $C \setminus Z$\uscolon


 
 \eenum
 in fact, for any extension of $B \subset G$ to a quasi-pinning\footnote{\label{foot:what-is-qs}We recall from \cite{SGA3IIInew}*{Expos\'{e}~XXIV, Sections 3.8--3.9} that a \emph{quasi-pinning} of a reductive group scheme $G$ is the datum of a Borel subgroup $B \subset G$, a maximal torus $T \subset B$, and, on the scheme $\underline{\Dyn}(G)$ of those $T$-roots that are simple with respect to $B$, a trivialization of the line bundle given by the universal root space that is simple with respect to $B$. A reductive group scheme is \emph{quasi-split} if it admits a quasi-pinning---over a semilocal base this amounts to having a Borel subgroup, but in general it is a more stringent condition.} of $G$, one may find the objects as above such that $\sB \subset \sG$ extends to a quasi-pinning of $\sG$ that  is compatible via $s$ with the quasi-pinning of $G$.
 
\epropt


\bpf
We decompose $\Spec(\cO)$ and $\Spec(R)$ into connected components to assume that $\cO$ and $R$ are domains, and then likewise assume that $A$ is a domain. If $A$ is of relative dimension $0$, then $R$ is a Dedekind domain, so, by \cite{Guo20}*{Theorem~1}, our torsor $E$ is trivial and we may choose $C  = \bA^1_R$ and $\sE \ce E_{\bA^1_R}$, the closed subscheme $Z$ being empty and $s$ being the zero section. 
 Thus, we may assume that $A$ is $\cO$-fiberwise of pure dimension $d > 0$.  Moreover, we localize $A$ and spread out to assume (abusively, from a notational standpoint) that $G$, $B$, and $E$ begin life over $A$ and, for the last aspect of the claim, that a fixed extension of $B \subset G$ to a quasi-pinning of $G$ likewise exists over $A$.



By \cite{SGA3IIInew}*{Expos\'{e}~XXVI, Corollaire~3.6, Lemme~3.20}, the quotient $E/B$ is representable by a projective 
$A$-scheme. Thus, due to the generic triviality of $E$ and the valuative criterion of properness, there is a closed subscheme $Y \subset \Spec(A)$ of codimension $\ge 2$ such that $(E/B)_{\Spec(A)\setminus Y}$ has a section that generically lifts to $E$, in other words, such that $E_{\Spec(A)\setminus Y}$ reduces to a generically trivial $B_{\Spec(A)\setminus Y}$-torsor $E^B$. Consider the torus 
\[
T \ce B/\sR_u(B) \qxq{and the induced $T_{\Spec(A)\setminus Y}$-torsor} E^T \ce E^B/\sR_u(B).
\]
Since $Y$ is of codimension $\ge 2$ in the regular scheme $\Spec(A)$, by \cite{CTS79}*{Corollaire~6.9}, there is a unique $T$-torsor $\wt{E^T}$ that extends $E^T$ to all of $\Spec(A)$. Since the Grothendieck--Serre conjecture is known for tori \cite{CTS87}*{Theorem~4.1~(i)} and $\wt{E^T}$ is generically trivial, the base change of $\wt{E^T}$ to $\Spec(R)$ is trivial. Thus, we may localize $A$ further around the maximal ideals of $R$ to assume that $\wt{E^T}$ is trivial, so that $E^T$ is also trivial and $E_{\Spec(A)\setminus Y}$ reduces to an $\sR_u(B)$-torsor.



We now apply \Cref{cor:main-structural} to obtain 
\bitem
\m
an affine open $U \subset \Spec(A)$ containing $\Spec(R)$; 

\m
an affine open $S \subset \bA^{d - 1}_\cO$; and

\m
a smooth $\cO$-morphism $U \ra S$ of pure relative dimension $1$ such that $Y \cap U$ is
 $S$-finite. 
\eitem
Since $R$ is a localization of the coordinate ring of $U$, we then set
\[
C \ce U \times_S \Spec(R) \qxq{and} Z \ce (Y \cap U) \times_S \Spec(R).
\]
The $R$-scheme $C$ comes equipped with an $R$-point $s \in C(R)$ induced by the diagonal of $\Spec(R)$ over $S$, and, by base change, \ref{CG-D}--\ref{CG-iv} hold.
Finally, we let $\sG$, $\sB$, and $\sE$ be the base changes to $C$ of the restrictions of $G$, $B$, and $E$ to $U$, so that, by construction, their $s$-pullbacks are $G$, $B$, and $E$, respectively. Since $U \subset \Spec(A)$ and, by construction, $E_{\Spec(A)\setminus Y}$ reduces to an $\sR_u(B)$-torsor, the restriction of $\sE$ to $C \setminus Z$ reduces to an $\sR_u(\sB)$-torsor. In particular, \ref{CG-v-new} and \ref{CG-v} hold. To likewise arrange the final claim about quasi-pinnings, simply take the fixed extension of $B \subset G$ to a quasi-pinning of $G$ and base change it to $C$.
\epf

\bremt
\label{rem:equichar}
Proposition \ref{cor:C-and-G} is significantly simpler in the case when $\cO$ is a field, in which the assumption that $G$ be quasi-split could be removed with $\sE$ in \ref{CG-v} being trivial over $C \setminus Z$. The point is that in this case, even if $Y$ is only of codimension $> 0$, as is immediate to arrange from the generic triviality of $E$, one nevertheless gets that $Y'$ in the proof of \Cref{cor:main-structural} is of codimension $> 1$ and the argument goes through without worrying about $B$ and with $\sE|_{C\setminus Z}$ even being trivial. 
\eremt


\section{Changing the relative curve to equate $\sG$ and $G_C$} \label{sec:fix-G}

To reach our main result on the Grothendieck--Serre conjecture we will gradually simplify the structure exhibited in \Cref{cor:C-and-G} and eventually reduce to studying $G$-torsors over $\bA^1_R$. As the first step towards this, in \Cref{cor:only-C} we reduce to the case when the $C$-group scheme $\sG$ appearing there is constant, that is, when $\sG$ is simply $G_C$. Although this step could be pursued more generally, here we simplify it by taking advantage of the quasi-splitness assumption. The following lemma is 
 close in spirit to \cite{OP01}*{Proposition~7.1}, \cite{PSV15}*{Proposition~5.1}, or \cite{Pan20b}*{Theorem~3.4}. 

\blemt \label{prop:finer-invariance}
For  a Noetherian semilocal ring $A$ whose localizations at prime ideals are geometrically unibranch, an ideal $I \subset A$, reductive $A$-groups $G$ and $G'$ that on geometric $A$-fibers have the same type, fixed quasi-pinnings of $G$ and $G'$ extending Borel $A$-subgroups $B \subset G$ and $B' \subset G'$, and an $A/I$-group isomorphism 
\[
\iota \colon G_{A/I} \isomto G'_{A/I} \qxq{respecting the quasi-pinnings, so that, in particular,} \iota(B_{A/I}) = B'_{A/I},
\]
there are
\benumr
\item
a faithfully flat, finite, \'{e}tale $A$-algebra $\wt{A}$ equipped with an $A/I$-point $a \colon \wt{A} \surjects A/I$\uscolon and

\item
an $\wt{A}$-group isomorphism $\wt{\iota} \colon G_{\wt{A}} \isomto G'_{\wt{A}}$ respecting the quasi-pinnings such that $a^*(\wt{\iota})=\iota$.
\eenum
\elemt


\bpf
We consider the functor $X$ that parametrizes those isomorphisms between the base changes of $G$ and $G'$ that preserve the (corresponding base changes of the) fixed quasi-pinnings. By \cite{SGA3IIInew}*{Expos\'{e}~XXIV, Corollaire~1.10, Section 3.10}, this $X$ is representable by a scheme that becomes constant \'{e}tale locally on $A$. Thus, by \cite{SGA3II}*{Expos\'{e}~X, Corollaire~5.14} (with \cite{EGAI}*{Corollaire~6.1.9}), the geometrically unibranch assumption ensures that the connected components of $X$ are open subschemes that are finite \'{e}tale over $A$. The $A/I$-point $\iota$ of $X$ meets finitely many such components, whose union is then the spectrum of a finite \'{e}tale $A$-algebra $\wt{A}$. By adjoining further components if needed, we may ensure that the closed $A$-fibers of $\wt{A}$ are nonzero, so that $\wt{A}$ is faithfully flat over $A$, as desired. 
\epf

\bpropt \label{cor:only-C}
For a semilocal Dedekind ring $\cO$, the semilocalization $R$ of a smooth $\cO$-algebra at finitely many primes $\fp$, 
a quasi-split reductive $R$-group $G$, a Borel $R$-subgroup $ B \subset G$, 
and a generically trivial $G$-torsor $E$,
there are
\benumr

\m 
a smooth, affine $R$-scheme $C$ of pure relative dimension $1$\uscolon


\m 
a section $s \in C(R)$\uscolon

\m 
an $R$-finite closed subscheme $Z \subset C$\uscolon


 \m 
a $G_C$-torsor $\sE$ whose $s$-pullback is $E$ such that $\sE$ reduces to an $\sR_u(B)$-torsor over $C \setminus Z$.
 \eenum
\epropt

\bpf
We decompose $\Spec(R)$ into connected components to assume that $R$ is a domain. By \Cref{cor:C-and-G}, there are such $C$, $s$, $Z$, and $\sE$, except that $\sE$ there is a torsor under a quasi-split reductive $C$-group scheme $\sG$ that may not be $G_C$ but that comes equipped with a Borel $C$-subgroup $\sB \subset \sG$ whose $s$-pullback is $B \subset G$ and that extends to a quasi-pinning of $\sG$ whose $s$-pullback is a fixed quasi-pinning of $G$ extending $B$. We replace $C$ by its connected component containing the image of $s$ to arrange that $C$ be connected. Thus, the geometric $C$-fibers of $\sG$ and $G_C$ are of constant types, so that, by the condition on the $s$-pullback, these types are the same. 

To replace $\sB \subset \sG$ by $B_C \subset G_C$, we first use prime avoidance \cite{SP}*{Lemma~\href{https://stacks.math.columbia.edu/tag/00DS}{00DS}} to construct the semilocalization $\Spec(A)$ of $C$ at the union of the closed points of $Z$ and of those of the image of $s$. Since $C$ is $R$-smooth, the ring $A$ is regular. The image of $s$ gives rise to a closed subscheme $\Spec(R) \subset \Spec(A)$ cut out by an ideal $I \subset A$ and, by assumption, $\sB_{A/I} \subset \sG_{A/I}$ agrees with $B \subset G$, and even the quasi-pinnings of $\sG$ and $G$ are likewise compatible. Thus, by \Cref{prop:finer-invariance}, there is a finite \'{e}tale $\Spec(A)$-scheme $\Spec(\wt{A})$ equipped with an $R$-point $\wt{s}$ lifting $s$ such that $\sB_{\wt{A}} \subset \sG_{\wt{A}}$ is isomorphic to $B_{\wt{A}} \subset G_{\wt{A}}$  compatibly with the fixed identification of $\wt{s}$-pullbacks. We may spread out 
\[
\Spec(\wt{A}) \ra \Spec(A) \qxq{to a finite \'{e}tale morphism} \wt{C} \ra C'
\]
for some affine open $C' \subset C$ that contains $Z$ and the image of $s$, while preserving an $\wt{s} \in \wt{C}(R)$ and an isomorphism between $\sB_{\wt{C}} \subset \sG_{\wt{C}}$ and $B_{\wt{C}} \subset G_{\wt{C}}$. 
To arrive at the desired conclusion, it then remains to replace $C$, $s$, $Z$, and $\sE$ by $\wt{C}$, $\wt{s}$, $Z \times_{C} \wt{C}$, and $\sE \times_C \wt{C}$, respectively. 
\epf





\section{A Lindel trick in the setting of Cohen--Macaulay relative curves}
\label{sec:Lindel}

Having corrected $\sG$, our next goal is to reduce to the case when the affine relative curve $C$ is $\bA^1_R$. 
For this, we need a suitable flat morphism $C \ra \bA^1_R$, whose construction is the goal of this section. We summarize the resulting relevant for us refinement of \Cref{cor:only-C} in \Cref{prop:better-C}.

To be able to later use patching to pass from $C$ to $\bA^1_R$, it is key to arrange that on some open subscheme $C' \subset C$ containing $Z$ our desired flat map $C \ra \bA^1_R$ does not change $Z$ in the sense that the latter is precisely the scheme-theoretic preimage in $C'$ of some closed subscheme $Z' \subset \bA^1_R$ to which $Z$ maps isomorphically. This is reminiscent of Lindel's insight \cite{Lin81}*{page~321, Lemma} that led to the resolution of the Bass--Quillen conjecture in the ``geometric'' case and says that an \'{e}tale map $B \ra A$ of local rings with the same residue field is an isomorphism modulo powers of a well-chosen element in the maximal ideal of $B$ (compare also with \cite{CTO92}*{Lemme~1.2} or \cite{CT95}*{Section 3.7 and proof of Theorem~5.1.1}). In our situation, however, there is a basic obstruction to the existence of $Z'$: if some residue fields of $R$ are finite, then $Z$ could have too many rational points to fit into $\bA^1_R$. The purpose of the following crucial statement essentially taken from the literature is to circumvent this obstacle.







\blemt \label{lem:ff-trick}
For a semilocal ring $R$, a quasi-projective, finitely presented $R$-scheme $C$, its $R$-finite closed subscheme $Z$, and an $s \in Z(R)$, there is a finite morphism $\wt{C} \ra C$ that is \'{e}tale at the points in $\wt{Z} \ce Z \times_{C} \wt{C}$ such that $s$ lifts to $\wt{s} \in \wt{C}(R)$ and, for every maximal ideal $\fm \subset R$, we have
\[
\#\{ z \in \wt{Z}_{k_\fm}\, \vert\, [k_z : k_\fm] = d \} < \#\{ z \in \bA^1_{k_\fm}\, \vert\, [k_z : k_\fm] = d \} \qxq{for every} d \ge 1
\]
\up{a vacuous condition if the residue field $k_\fm$ is infinite}.
\elemt

Lemmas like \ref{lem:ff-trick} and \ref{lem:Lindel-variant} below go back to Panin's \cite{Pan15}*{Lemma~14.1 and Theorem~11.7}.

\bpf
The lemma is a variant of, for instance, \cite{Pan19a}*{Lemma~6.1} or \cite{Fed21b}*{Lemma~4.5}, and we will prove it by using similar arguments as there due to Panin. Since $R$ is semilocal, the finite $R$-scheme $Z$ has finitely many closed points, which all lie over maximal ideals of $R$. Thus, we begin by using \Cref{lem:avoidance} to construct the semilocalization $S$ of $C$ at the closed points of $Z$, so that $Z$ is also a closed subscheme of $S$ and $s \in S(R)$. It then suffices to construct a finite \'{e}tale $S$-scheme $\wt{S}$ such that $s$ lifts to an $R$-point $\wt{s} \in \wt{S}(R)$ and the preimage $\wt{Z} \subset \wt{S}$ of $Z$ satisfies the displayed inequalities: indeed, once this is done, we may first spread $\wt{S}$ out  to a finite \'{e}tale scheme over an open neighborhood of $S$ in $C$ and then form its schematic image \cite{SP}*{Lemma~\href{https://stacks.math.columbia.edu/tag/01R8}{01R8}} in the factorization supplied by the Zariski Main Theorem \cite{EGAIV4}*{Corollaire~18.12.13} to further extend to a desired finite morphism $\wt{C} \ra C$. 

We view $s$ as a closed subscheme $\Spec(R) \subset Z$ and we list the closed points of $Z$ (that is, of $S$):
\bitem
\m
the closed points $y_1, \dotsc, y_{m}$ of $Z$ not in $s$ with an infinite residue field;

\m
the closed points $z_1, \dotsc, z_{n}$ of $Z$ not in $s$ with a finite residue field;

\m
the closed points $y_1', \dotsc, y'_{m'}$ of $s$ with an infinite residue field;

\m
the closed points $z_1', \dotsc, z'_{n'}$ of $s$ with a finite residue field.

\eitem
For any $N > 1$, we may choose monic polynomials
\bitem
\m
$f_{y_i} \in k_{y_i}[t]$ that are products of $N$ distinct linear factors; and 

\m
$f_{z_j} \in k_{z_j}[t]$ that are irreducible of degree $N$. 
\eitem
Likewise, we may choose a monic polynomial $f_{s} \in tR[t]$ of degree $N$ such that 
\bitem
\m
the image of $f_{s}$ in each $k_{y_i'}[t]$ is a product of $N$ distinct linear factors; and

\m
the image of $f_{s}$ in each $k_{z_j'}[t]$ is a product of $t$ and an irreducible polynomial not equal to $t$. 
\eitem
Finally, since $s \sqcup \bigsqcup_{i = 1}^{m} y_i \sqcup \bigsqcup_{j = 1}^{n} z_j$ is a closed subscheme of $S$, by lifting coefficients we may choose a monic polynomial $f \in \Gamma(S, \sO_S)[t]$ of degree $N$ that restricts to $f_{y_i}$ on each $y_i$, to $f_{z_j}$ on each $z_j$, and to $f_{s}$ on $s$. This $f$ defines a finite \'{e}tale $S$-scheme $\wt{S}$, which, by construction, is equipped with an $R$-point $\wt{s} \in \wt{S}(R)$ lifting $s$ (cut out by the factor $t$ of $f_{s}$) and is such that  the number of closed points with finite residue fields in the preimage $\wt{Z} \subset \wt{S}$ of $Z$ stays bounded as $N$ grows but, except for the points in $\wt{s}$, the cardinalities of these residue fields grow uniformly. Thus, since, for a finite field $\bF$, the number of closed points of $\bA^1_\bF$ with a given residue field grows unboundedly together with the degree of that residue field over $\bF$, for large $N$ our $\wt{S}$ meets the requirements. 
\epf

We turn to the Lindel trick in our setting, namely, to building the desired flat map $C \ra \bA^1_R$ in \Cref{lem:Lindel-variant}. Its numerous variants appeared in works of Panin, for instance, in \cite{OP99}*{Section 5}, \cite{PSV15}*{Theorem~3.4}, or \cite{Pan19a}*{Theorem~3.8}, but with the more stringent smoothness assumption on $C$, and  preparation lemmas of similar flavor can be traced back at least to \cite{Gab94b}*{Lemma~3.1} (compare also with \cite{CTHK97}*{Theorem~3.1.1}). 
As we show, Cohen--Macaulayness of $C$ suffices. 
The argument uses the following simple lemma that characterizes residue fields of closed points on smooth curves. 

\blemt \label{lem:curve-pts}
For a field $k$, a smooth connected $k$-curve $C$, and a closed point $c \in C$, the  extension $k_c/k$ is generated by a single element, that is, $k_c$ is the residue field of a closed point of $\bA^1_k$.  
\elemt

\bpf
By \cite{EGAIV4}*{Corollaire~17.11.4}, an open neighborhood $U \subset C$ of $c$ has an \'{e}tale $k$-morphism 
\[
U \ra \bA^1_k.
\]
Thus, there is a subextension $\ell/k$ of $k_c/k$ generated by a single element with $k_c/\ell$ separable. By the primitive element theorem, we need to check that this forces $k_c/k$ to only have finitely many subextensions $k'/k$. Since there are finitely many possibilities for $k' \cap \ell$, we replace $k$ by $k' \cap \ell$ to reduce to considering those $k'$ for which $k' \cap \ell = k$. Like any finite separable extension, the separable closure of $k$ in $k_c$ has only finitely many subextensions. Thus, there are finitely many possibilities for the maximal separable subextension $k''/k$ of $k'/k$. By replacing $k$ by $k''$ and $\ell$ by $k''\ell$, we therefore reduce to the case when $k'/k$ is purely inseparable. Then the subextension $k'\ell/\ell$ of the separable extension $k_c/\ell$ is also purely inseparable, to the effect that $k' \subset \ell$. However, $\ell/k$ is generated by a single element, so, by the primitive element theorem, it has only finitely many subextensions.
\epf

\blemt \label{lem:Lindel-variant}
For 
\bitem
\m
a  semilocal ring $R$\uscolon

\m
a flat, affine $R$-scheme $C$ with Cohen--Macaulay fibers of pure dimension $1$\uscolon 

\m
$R$-finite closed subschemes $Y \subset C$ and $Z \subset C^\sm$ such that, for every maximal ideal $\fm \subset R$,
\[
\qq \#\{ z \in Z_{k_\fm}\, \vert\, [k_z : k_\fm] = d \} < \#\{ z \in \bA^1_{k_\fm}\, \vert\, [k_z : k_\fm] = d \} \qxq{for every} d \ge 1
\]
\up{a vacuous condition if the residue field $k_\fm$ is infinite}\footnote{If every closed point of $Y$ lies on $Z$, then it suffices to require the same with nonstrict inequalities $\le$ instead.}\uscolon 
\eitem
there are 
\benumr
\m \label{LV-i}
an affine open $C' \subset C$ containing $Y$ and $Z$\uscolon

\m\label{LV-ii}
a quasi-finite, flat $R$-map $C' \ra \bA^1_R$ that maps $Z$ isomorphically onto a closed subscheme 
\[
Z' \subset \bA^1_R \qxq{
such that} Z \cong Z' \times_{\bA^1_R} C'; 
\]
\eenum
so that, in particular, $C' \ra \bA^1_R$ is \'{e}tale along $Z$ and, for every $n \ge 0$,  maps the $n$-th infinitesimal neighborhood of $Z$ in $C'$ isomorphically to the $n$-th infinitesimal neighborhood of $Z'$ in $\bA^1_R$.
\elemt

\emph{Proof.}\footnote{We loosely follow \cite{Pan19a}*{proof of Theorem~3.8}, with several improvements and simplifications whose purpose is to avoid assuming that $C$ be $R$-smooth or that $R$ be the semilocal ring at finitely many closed points of a smooth variety over a field. Notably, in \Cref{rem:CH-trick} we give a more direct and more general argument for the final portion of \emph{loc.~cit.} 
}
\addtocounter{footnote}{-5}
\renewcommand{\thefootnote}{\fnsymbol{footnote}}
The \'{e}taleness follows from the flatness and the isomorphy over $Z'$ of the map $C' \ra \bA^1_R$, and it implies the infinitesimal neighborhood aspect.
 For the rest, 
every closed point $z \in Z$ lies over some maximal ideal $\fm \subset R$ and, since $z \in C^\sm_{k_\fm}$, the ideal sheaf $\sI_z \subset \sO_{C_{k_\fm}}$ is generated at $z$ by a uniformizer $u_z \in \sO_{C_{k_\fm},\, z}$. Consequently, by \cite{BouAC}*{Chapitre IX, Section 3, num\'{e}ro~3, Th\'{e}or\`{e}me~1}, the thickening\footnote[2]{\emph{Added after publication.} As Gabber pointed out, the isomorphism $\eps_z \simeq \Spec(k_z[u_z]/(u_z^2))$ deduced here from the Cohen structure theorem is \emph{a priori} not an isomorphism of $k_\fm$-schemes unless the field $k_\fm$ is perfect (in which case $k_z/k_\fm$ is a finite separable extension and the claimed $k_\fm$-scheme isomorphism also follows from \cite{EGAIV4}*{Proposition~17.5.3} combined with the insensitivity of the \'{e}tale site to nilpotents). This does not affect the subsequent argument because, for imperfect $k_\fm$, we may nevertheless construct the closed immersions $\eps_z \hra \bA^1_{k_\fm}$ comprising $j$ by applying the geometric presentation theorem \cite{CTHK97}*{Theorem 3.1.1} to each $\eps_z \subset C_{k_\fm}^\sm$ and postcomposing with suitable changes of coordinates $t \mapsto t + \gA$ of $\bA^1_{k_\fm}$ for $\gA \in k_\fm$ to make the resulting closed immersions $\eps_z \hra \bA^1_{k_\fm}$ have disjoint images.}
\[
\eps_z \ce \underline{\Spec}_{\sO_{C_{k_\fm}}}(\sO_{C_{k_\fm}}/\sI_z^2) \qxq{is isomorphic to} \Spec(k_z[u_z]/(u_z^2)).
\]
Letting $y$ range over the closed points of $Y$ not in $Z$ and $z$ range over the closed points of $Z$, we set
\[
\tst \eps_Y \ce \bigsqcup_y y \subset C, \q \eps_Z \ce \bigsqcup_z \eps_z \subset C, \qxq{and} \eps \ce \eps_Y \sqcup \eps_Z = \bigsqcup_y y \sqcup \bigsqcup_z \eps_z \subset C.
\]
By \Cref{lem:curve-pts} and the assumption on the numbers of points of $Z_{k_\fm}$, we may find an $R$-morphism
\[
\tst  j\colon \eps \ra \bigsqcup_{\fm 
} \bA^1_{k_\fm} \subset \bA^1_R \qxq{that restricts to a closed immersion} \eps_Z \hra \bA^1_R
\] 
and, for each $\fm$, maps the points of $\eps_Y$ above $\fm$ to an $k_\fm$-point of $\bA^1_{k_\fm} \setminus \eps_Z$. We fix 
 two disjoint sets of closed points
\[
d_1, \dotsc, d_n \in C\setminus(Y \cup Z) \qxq{and} d'_1, \dotsc, d'_{n'} \in C \setminus(Y \cup Z)
\]
lying over maximal ideals of $R$ such that each set jointly meets every irreducible component of every closed $R$-fiber of $C$. Since $\eps \cup \bigsqcup_{i = 1}^n d_i \cup \bigsqcup_{i = 1}^{n'} d_i'$ is a closed subscheme of $C$,
we may find 
an $s \in \Gamma(C, \sO_C)$ that 
\bitem
\m
vanishes on every $d_i$ but does not vanish on any $d_i'$;

\m
 on $\eps$ equals the $j$-pullback of the coordinate of $\bA^1_R$. 
\eitem
By mapping the coordinate of $\bA^1_R$ to $s$, we obtain an $R$-morphism 
\[
\pi \colon C \ra \bA^1_R.
\]
The behavior of $s$ at $d_i$ and $d_i'$ ensures that the locus where $\pi$ is quasi-finite, which, by \cite{SP}*{Lemma~\href{https://stacks.math.columbia.edu/tag/01TI}{01TI}}, is an open of $C$, contains every closed $R$-fiber of $C$. In particular, we may use prime avoidance \cite{SP}*{Lemma~\href{https://stacks.math.columbia.edu/tag/00DS}{00DS}} to replace $C$ by some affine open subset containing $Y$ and $Z$ (equivalently, containing the closed points of $Y$ and $Z$) to arrange that $\pi$ is quasi-finite. 

Since $C$ is $R$-flat with Cohen--Macaulay fibers of pure dimension $1$, the flatness criteria \cite{EGAIV2}*{Proposition~6.1.5}, \cite{EGAIV3}*{Corollaire~11.3.11} ensure that $\pi$ is  flat. By construction $\pi|_{\eps} = j$, so, by checking on the closed $R$-fibers, \cite{EGAIV4}*{Th\'{e}or\`{e}me~17.11.1} shows that $\pi$ is \'{e}tale around $Z$. Since $Z_{k_\fm}$ and $\eps_Z$ have the same underlying reduced subscheme $\bigsqcup_z z$, the agreement with $j$ also shows that $\pi|_{Z_{k_\fm}}$ is a closed immersion. 
Since $Z$ is $R$-finite, Nakayama lemma \cite{SP}*{Lemma~\href{https://stacks.math.columbia.edu/tag/00DV}{00DV}} then ensures that $\pi|_Z$ is also a closed immersion, so that $\pi$ maps $Z$ isomorphically onto a closed subscheme $Z' \subset \bA^1_R$. 

A section of a separated, \'{e}tale morphism is an isomorphism onto a clopen subscheme, so the \'{e}taleness of $\pi$ around $Z$ gives a decomposition 
\[
\pi\i(Z') = Z \sqcup Z''
\]
for some $R$-quasi-finite closed subscheme $Z'' \subset C$. By the agreement with $j$, the image under $\pi$ of every closed point of $Y$ not in $Z$ does not lie in $Z'$, so that $Y \cap Z'' = \emptyset$. Thus, prime avoidance \cite{SP}*{Lemma~\href{https://stacks.math.columbia.edu/tag/00DS}{00DS}} supplies a global section of $C$ that vanishes on $Z''$ but does not vanish at any closed point of $Y$ or $Z$. By inverting this section, we obtain the desired affine open $C' \subset C$. 
\QED

\renewcommand{\thefootnote}{\arabic{footnote}}
\addtocounter{footnote}{5}

\bremt \label{rem:CH-trick}
If $\Spec(R)$ is connected,  
then any $R$-(finite locally free) closed subscheme $Z' \subset \bA^1_R$ is cut out by a monic polynomial. This holds for any ring $R$ with a connected spectrum: the coordinate $t$ of $\bA^1_R$ acts by multiplication on the projective $R$-module $\Gamma(Z', \sO_{Z'})$, the characteristic polynomial of this action is monic and cuts out an $R$-(finite locally free) closed subscheme $H \subset \bA^1_R$, and Cayley--Hamilton implies that $Z' \subset H$ inside $\bA^1_R$, so, by comparing ranks over $R$, even $Z' = H$.
\eremt

We now refine \Cref{cor:only-C} to the following statement adapted to passing to $\bA^1_R$ via patching. 

\bpropt \label{prop:better-C}
For a semilocal Dedekind ring $\cO$, the localization $R$ of a smooth $\cO$-algebra at finitely many primes $\fp$, 
a quasi-split reductive $R$-group $G$, a Borel $R$-subgroup $ B \subset G$, 
and a generically trivial $G$-torsor $E$,
there are
\benumr
\m \label{bC-i}
a smooth, affine $R$-scheme $C$ of pure relative dimension~$1$\uscolon


\m \label{bC-ii}
a section $s \in C(R)$\uscolon

\m \label{bC-iii}
an $R$-finite closed subscheme $Z \subset C$\uscolon

 \m \label{bC-iv}
a $G_C$-torsor $\sE$ whose $s$-pullback is $E$ such that $\sE$ reduces to a $\sR_u(B)$-torsor over $C \setminus Z$\uscolon

\m \label{bC-vi}
a flat $R$-map $C \ra \bA^1_R$ that maps $Z$ isomorphically onto a closed subscheme $Z' \subset \bA^1_R$ with
\[
\qq Z \cong Z' \times_{\bA^1_R} C.
\]
 \eenum	
\epropt

\bpf
\Cref{cor:only-C} supplies $C$, $s$, $Z$, and $\sE$ that satisfy the present \ref{bC-i}--\ref{bC-iv}. 
We view $s$ as a closed subscheme of $C$ and we apply \Cref{lem:ff-trick} to the $R$-finite closed subscheme $(Z \cup s)^\red$ of $C$ to see that we may change $C$ to assume that, in addition, for every maximal ideal $\fm \subset R$, 
\[
\#\{ z \in  Z_{k_\fm} \, \vert\, [k_z : k_\fm] = d \} < \#\{ z \in \bA^1_{k_\fm}\, \vert\, [k_z : k_\fm] = d \} \qxq{for every} d \ge 1.
\]
This allows us to apply \Cref{lem:Lindel-variant} with $Y = s$ to shrink $C$ and to arrange \ref{bC-vi}.
\epf


\section{Descending to $\bA^1_R$ via patching} \label{sec:to-A1}

With the suitable flat map $C \ra \bA^1_R$ already built in \Cref{prop:better-C}, 
descending the $G_C$-torsor $\sE_C$ to $\bA^1_R$ concerns patching along the closed subscheme $Z$. Since our $Z$ need not be cut out by a single equation (relatedly, $C\setminus Z$ need not be affine), this patching is slightly more delicate than its most frequently encountered instances. Its precise statement is captured by the following lemma, which follows from more general results of Moret-Bailly \cite{MB96} (for our purposes, we could also get by with the more basic patching of Ferrand--Raynaud \cite{FR70}*{Proposition~4.2}).

\blemt \label{lem:LMB-patching}
Let $S' \ra S$ be an affine, flat scheme map whose base change to a closed subscheme $Z \subset S$ cut out by a quasi-coherent ideal sheaf of finite type is an \emph{isomorphism} and let $U' \ra U$ be the base change to  $U \ce S\setminus Z$. For a quasi-affine, flat, finitely presented $S$-group scheme $G$, base change induces an equivalence from the category of $G$-torsors to the category of triples consisting of a $G_{S'}$-torsor, a $G_{U}$-torsor, and a $G_{U'}$-torsor isomorphism between the two base changes to $U'$. 
\elemt

Of course, the isomorphism condition $Z \times_S S' \isomto Z$ ensures that $S'$ and $U$ jointly cover $S$.

\bpf
By \cite{SP}*{Theorem \href{https://stacks.math.columbia.edu/tag/06FI}{06FI}}, the classifying $S$-stack $\bbB G$ is algebraic and, by descent, its diagonal inherits quasi-affineness from $G$. Thus, the claim is a special case of \cite{MB96}*{Corollaire~6.5.1~(a)}.
\epf

To be able to apply \Cref{lem:LMB-patching} in our setting, we need to descend $\sE_{U \setminus Z}$ to a $G$-torsor over $\bA^1_R \setminus Z'$. To achieve this, we will use the following excision result that is similar (but simpler) than its counterparts that recently appeared in \cite{flat-purity}*{Theorem~5.4.4} and in \cite{Hitchin-torsors}*{Section 2.3}.

\blemt \label{lem:excision}
Let $S' \ra S$ be a flat morphism of affine, Noetherian schemes whose base change to a closed subscheme $Z \subset S$ is an isomorphism, and let $U' \ra U$ be the base change to $U \ce S \setminus Z$. 
\benum
\m \label{E-a}
For a quasi-coherent $\sO_S$-module $\sF$ \up{or even a complex of such $\sO_S$-modules}, we have 
\[
\q R\Gamma_Z(S, \sF) \isomto R\Gamma_{Z}(S', \sF_{S'}).
\] 

\m \label{E-b}
For an affine, smooth $S$-group \up{resp.,~$U$-group} $F$ with a filtration 
\[
\qq F = F_0 \supset F_1 \supset \dotsc \supset F_n = 0
\]
by normal, affine, smooth $S$-subgroups \up{resp.,~$U$-subgroups} such that, for all $i \ge 0$, the quotient $F_i/F_{i + 1}$ is a vector group associated to a vector bundle on $S$ \up{resp.,~such that the vector group $F_i/F_{i + 1}$ is also central in $F/F_{i + 1}$}, the map
\[
\q H^1(U, F) \ra H^1(U', F) \qx{has trivial kernel \up{resp.,~is surjective}.}
\]
\eenum
\elemt

\bpf \hfill
\benum
\m
We let $A$ and $A'$ be the coordinate rings of $S$ and $S'$, respectively. By \cite{SP}*{Lemmas \href{https://stacks.math.columbia.edu/tag/0ALZ}{0ALZ} and \href{https://stacks.math.columbia.edu/tag/0955}{0955}},
\[
\qq R\Gamma_Z(S, \sF) \tensor^\bL_A A' \isomto R\Gamma_{Z}(S', \sF_{S'}).
\]
Thus, since $A'$ is $A$-flat, to obtain \ref{E-a} it remains to note that, by \cite{SP}*{Lemma~\href{https://stacks.math.columbia.edu/tag/05E9}{05E9}}, we have
\[
\qq H^i_Z(S, \sF) \isomto H^i_Z(S, \sF) \tensor_A A' \qxq{for all} i \in \bZ. 
\]

\m
In the case when $F$ is an $S$-group, the vanishing of quasi-coherent cohomology of affine schemes and the assumed filtration show that both $H^1(S, F)$ and $H^1(S', F)$ vanish. Thus, the assertion about the kernel simply amounts to the claim that every $F_U$-torsor that trivializes over $U'$ extends to an $F$-torsor. This, however, is immediate from \Cref{lem:LMB-patching}. 

For the surjectivity assertion, we will induct on $n$. We begin with the case $n = 1$, in which $F$ itself is the vector group associated to some vector bundle $\sF$ on $U$. By applying \ref{E-a} to $j_*(\sF)$, where $j \colon U \hra S$ is the indicated open immersion, and again using the vanishing of quasi-coherent cohomology of affine schemes, we find that, for all $i\ge 1$, even
\[
\xymatrix{
\qq\ \ H^i(U, F) \cong H^i(U, \sF) \cong H^{i + 1}_{Z}(S, j_*(\sF)) \ar[r]_-{\sim}^-{\x{\ref{E-a}}} & H^{i + 1}_{Z}(S', j_*(\sF)) \cong H^i(U', \sF_{U'}) \cong  H^i(U', F).
}
\]
For the inductive step, we assume that $n > 1$ and combine the inductive hypothesis, the preceding display for $F_{n - 1}$, and the nonabelian cohomology sequences \cite{Gir71}*{Chapitre~IV, Remarque~4.2.10} of a central extension to obtain the following commutative diagram with exact rows:
\[
\qq\xymatrix{
 H^1(U, F_{n-1}) \ar[d]^{\sim} \ar[r] & H^1(U, F) \ar[r] \ar[d] & H^1(U, F/F_{n-1}) \ar@{->>}[d] \ar[r] & H^2(U, F_{n -1}) \ar[d]^{\sim} \\
 H^1(U', F_{n-1})  \ar[r] & H^1(U', F) \ar[r]  & H^1(U', F/F_{n-1}) \ar[r] & H^2(U', F_{n - 1}).
}
\]
We fix an $\gA' \in H^1(U', F)$ that we wish to lift to $H^1(U, F)$ and note that, by a diagram chase, there at least is an $\gA \in H^1(U, F)$ whose image in $H^1(U', F/F_{n - 1})$ agrees with that of $\gA'$. Every inner fpqc form of $F$ inherits an analogous filtration, even with the same subquotients $F_i/F_{i + 1}$, so the change of origin bijections \cite{Gir71}*{Chapitre~III, Proposition~2.6.1~(i)} allow us to twist $F$ and reduce to the case when the common image of $\gA$ and $\gA'$ in $H^1(U', F/F_{n - 1})$ vanishes. In this case, however, the surjectivity of the left vertical arrow suffices.
 \qedhere
\eenum
\epf

\begt \label{eg:parabolic}
For example, $F$ in \Cref{lem:excision}~\ref{E-b} could be the unipotent radical $\sR_u(P)$ of a parabolic $S$-subgroup (resp.,~$U$-subgroup) $P$ of a reductive $S$-group (resp.,~$U$-group) $G$: in this case, \cite{SGA3IIInew}*{Expos\'{e}~XXVI, Proposition~2.1} supplies the desired filtration.
\eegt

We can now reduce to the case when  the relative curve $C$ in \Cref{prop:better-C} is $\bA^1_R$.

\bpropt \label{cor:only-A1}
For a semilocal Dedekind ring $\cO$, the semilocalization $R$ of a smooth $\cO$-algebra at finitely many primes $\fp$, 
a quasi-split reductive $R$-group $G$, and a generically trivial $G$-torsor $E$, there are 
\benumr
\m
a closed subscheme $Z \subset \bA^1_R$ that is finite over $R$\uscolon

\m
 a $G_{\bA^1_R}$-torsor $\sE$ whose pullback along the zero section is $E$ such that $\sE$ is trivial over $\bA^1_R \setminus Z$.
 \eenum
\epropt

\bpf
Let $B \subset G$ be a Borel $R$-subgroup. \Cref{prop:better-C} supplies a quasi-finite, affine, flat $R$-morphism $\pi \colon C \ra \bA^1_R$ whose base change to an $R$-finite closed subscheme $Z \subset \bA^1_R$ (called $Z'$ there) is an isomorphism, as well as an $s \in C(R)$ and a $G_C$-torsor $\wt{\sE}$ (called $\sE$ there) with $s$-pullback $E$ such that $\wt{\sE}$ reduces to a $\sR_u(B)$-torsor over $C \setminus \pi\i(Z)$.  By \Cref{lem:excision}~\ref{E-b} and \Cref{eg:parabolic}, this $\sR_u(B)$-torsor descends to a $\sR_u(B)$-torsor over $\bA^1_R \setminus Z$, so $\wt{\sE}_{C \setminus \pi\i(Z)}$ descends to a $G_{\bA^1_R \setminus Z}$-torsor. The patching lemma \ref{lem:LMB-patching} then ensures that $\wt{\sE}$ itself descends to a $G_{\bA^1_R}$-torsor $\sE$ that reduces to a $\sR_u(B)$-torsor over $\bA^1_R \setminus Z$. By postcomposing  with a change of coordinate automorphism of $\bA^1_R$ to ensure that $s$ map to the zero section of $\bA^1_R$,  we make the pullback of $\sE$ along the zero section be $E$. Finally, we apply \Cref{lem:avoidance} to $\bP^1_R$ to enlarge $Z \subset \bA^1_R$ to be a hypersurface in $\bP^1_R$. This ensures that $\bA^1_R \setminus Z$ is affine, so that, due to the filtration of $\sR_u(B)$ by vector groups as in	 \Cref{eg:parabolic} and the vanishing of quasi-coherent cohomology of affine schemes, our  $\sE_{\bA^1_R \setminus Z}$ is even trivial.
\epf











\section{The analysis of torsors over $\bA^1_R$} \label{sec:A1-analysis}

Our final task is to study generically trivial torsors over $\bA^1_R$, which may be viewed as a problem of Bass--Quillen type beyond the case of vector bundles (that is, beyond $\GL_n$-torsors). Although one may carry it out much more generally, this analysis is significantly simpler for semisimple, simply-connected, totally isotropic reductive groups $G$. Total isotropicity, defined in \Cref{def:tot-iso}, is not quite unexpected, for instance, it also seems necessary for generalizing  the Bass--Quillen conjecture to $G$-torsors, see \cite{AHW18}*{Theorem~3.3.7} and \cite{AHW20}*{Theorem~2.4} for such results when $G$ comes from a ground field.

\bdt \label{def:tot-iso}
A reductive group scheme $G$ over a scheme $S$ is \emph{totally isotropic} if for every $s \in S$ 
and the canonical decomposition 
\[
\tst G^\ad_{\sO_{S,\, s}} \cong \prod_i G_i
\]
from \cite{SGA3IIInew}*{Expos\'{e}~XXIV, Proposition~5.10~(i)} of $G^\ad_{\sO_{S,\,s}}$ into the product of Weil restrictions $G_i$ of adjoint groups with simple geometric fibers over connected finite \'{e}tale covers of $\sO_{S,\,s}$, each $G_i$ is isotropic in the sense that it contains a nontrivial split torus $\bG_{m,\, \sO_{S,\,s}}$ (equivalently, each $G_i$ has a proper parabolic $\sO_{S,\,s}$-subgroup, see \cite{SGA3IIInew}*{Expos\'{e}~XXVI, Corollaire~6.12}). 
\edt

\begt\label{eg:qsplit}
The parenthetical reformulation implies that every quasi-split $G$ is totally isotropic.
\eegt

By the following lemma, which generalizes the main result of \cite{Tsy19}, for analyzing a torsor over $\bA^1_R$ the key is to extend it 
to a torsor over $\bP^1_R$ in such a way that the latter be trivial on closed $R$-fibers.

\blemt \label{lem:Fedorov}
For a semilocal ring $R$ and a reductive $R$-group $G$ such that $\rad(G)$ splits over a finite \'{e}tale cover of $\Spec(R)$ \up{a vacuous condition if $R$ is normal, or if $G$ is split, or if $\rad(G)$ is of rank $\le 1$, for instance, if $G$ is semisimple}, every $G_{\bP^1_R}$-torsor $\sE$ whose base change to $\bP^1_{k_\fm}$ is trivial for every maximal ideal $\fm \subset R$ is the base change of a $G$-torsor. 
\elemt

\bpf
By \cite{Gil21}*{Theorem 1.1 and Corollary 4.3}, the assumption on $\rad(G)$ is equivalent to requiring that $G$ may be embedded as a closed subgroup of some $\GL_{n,\, R}$ and it holds in the indicated parenthetical cases. 
For the rest, we first use a limit argument to reduce to Noetherian $R$ and then pass to connected components to also assume that $\Spec(R)$ is connected. Moreover, we begin with the case $G = \GL_{n,\, R}$, in which we may regard $\sE$ as a vector bundle of rank $n$. 

In this vector bundle case, 
\[
\sV \ce \sH om_{\sO_{\bP^1_R}}(\sO_{\bP^1_R}^{\oplus n}, \sE) \cong \sE^{\oplus n}
\]
is also a vector bundle on $\bP^1_R$. By \cite{EGAIII1}*{Th\'{e}or\`{e}me~3.2.1}, the $R$-module $V \ce \Gamma(\bP^1_R, \sV)$ is finite. By assumption, $\sE|_{\bP^1_{k_\fm}}$ is trivial for every maximal ideal $\fm \subset R$, so for such $\fm$ we choose an isomorphism 
\[
\sO_{\bP^1_{k_\fm}}^{\oplus n} \isomto \sE|_{\bP^1_{k_\fm}}, \qxq{which corresponds to some} v_\fm \in \Gamma(\bP^1_{k_\fm}, \sV|_{\bP^1_{k_\fm}}).
\]
Likewise, each $\sV|_{\bP^1_{k_\fm}}$ is trivial, so $H^1(\bP^1_{k_\fm}, \sV|_{\bP^1_{k_\fm}}) \cong 0$. Thus, by cohomology and base change \cite{EGAIII1}*{Proposition~4.6.1}, there is a $\wt{v}_\fm \in V/\fm V$ that maps to $v_\fm$. Since $R$ is semilocal and $\fm$ ranges over its maximal ideals, we may then find a $v \in V$ that maps to all the $\wt{v}_\fm$, so also to all the $v_\fm$. By construction and the Nakayama lemma \cite{SP}*{Lemma~\href{https://stacks.math.columbia.edu/tag/00DV}{00DV}}, the  $\sO_{\bP^1_R}$-module homomorphism 
\[
\sO_{\bP^1_R}^{\oplus n} \ra \sE
\]
 corresponding to $v$ is surjective at every closed point, so it is surjective. Cayley--Hamilton \cite{SP}*{Lemma~\href{https://stacks.math.columbia.edu/tag/05G8}{05G8}} then ensures that this homomorphism is an isomorphism, so $\sE$ is trivial, as desired.

To deduce the general case, we use our closed embedding $G \hra \GL_{n,\, R}$. Namely, the settled case of $\GL_{n,\, R}$ and the nonabelian cohomology sequence \cite{Gir71}*{Chapitre~III, Proposition~3.2.2} show that our $G_{\bP^1_R}$-torsor $\sE$ comes from a some $\bP^1_R$-point of $\GL_{n,\, R}/G$. However, $G$ is reductive, so, by \cite{Alp14}*{Theorem~9.4.1 and Theorem~9.7.5}, this quotient is affine. By \cite{MFK94}*{Proposition~6.1} (to reduce to an $R$-fiber), this means that the only $R$-morphisms from $\bP^1_R$ to $\GL_{n,\, R}/G$ are constant, in particular, that our $\bP^1_R$-point comes from an $R$-point. This then implies that our $G_{\bP^1_R}$-torsor $\sE$ is the base change of a $G$-torsor, as desired. 
\epf

The preceding lemma leads to the triviality of generically trivial reductive group torsors over $\bA^1_R$ as follows. 
Examples from \cite{Fed16a} show that without some isotropicity condition on $G$ such triviality does not hold.

\bpropt \label{prop:A1-trivial}
For a semilocal ring $R$ and a totally isotropic, semisimple, simply connected reductive $R$-group $G$,
every $G_{\bA^1_R}$-torsor $\sE$ that is trivial away from an $R$-finite closed subscheme $Z \subset \bA^1_R$ is trivial.


\epropt


\bpf
Due to the canonical decomposition \cite{SGA3IIInew}*{Expos\'{e}~XXIV, Sections 5.2--5.3, Proposition~5.10~(i)}, we may assume that 
\[
G \cong \Res_{R'/R} G'
\]
for a finite \'{e}tale $R$-algebra $R'$ and a semisimple, simply-connected $R'$-group $G'$ whose geometric $R'$-fibers are simple (in the sense that the Dynkin diagrams of these geometric fibers are connected). By \cite{SGA3IIInew}*{Expos\'{e}~XXIV, Proposition~8.4}, we may then replace $R$ and $G$ by $R'$ and $G'$, respectively, to assume that the geometric $R$-fibers of $G$ are simple.

We let $t$ be the inverse of the coordinate on $\bA^1_R$ and consider $R\llb t\rrb$ as the completion of $\bP^1_R$ along infinity. Due to its $R$-finiteness, $Z$ is closed in $\bP^1_R$, so its pullback to $\Spec(R\llb t \rrb)$ is also closed and hence is even empty because it does not meet the  locus $\{ t = 0\}$. Thus, we may use formal glueing supplied by, for instance, \cite{Hitchin-torsors}*{Lemma~2.2.11~(b)} (or \Cref{lem:LMB-patching} when $R$ is Noetherian) to extend $\sE$ to a $G_{\bP^1_R}$-torsor $\ov{\sE}$ by glueing $\sE$ with the trivial $G_{R\llb t \rrb}$-torsor. It suffices to argue that we can glue like this so that $\ov{\sE}_{\bP^1_{k_\fm}}$ be trivial for every maximal ideal $\fm \subset R$: \Cref{lem:Fedorov} will then imply that $\ov{\sE}$ is the base change of its pullback by the section at infinity, and hence that $\ov{\sE}$ and $\sE$ are trivial. 

Explicitly,  the elements of $G(R\llp t \rrp)/G(R\llb t \rrb)$ give rise to all the possible glueings of $\sE$ and the trivial $G_{R\llb t\rrb}$-torsor to a $G_{\bP^1_R}$-torsor, and likewise over the residue fields $k_\fm$ of $R$. We will first build a trivial $G_{\bP^1_{k_\fm}}$-bundle $\ov{\sE}_{\bP^1_{k_\fm}}$ from $\sE_{\bP^1_{k_\fm}}$ by such a glueing for every maximal ideal $\fm \subset R$ and then argue that these glueings come from a single glueing over $R$. These steps reduce, respectively, to the following:
\benuma
\m \label{Gille-1}
for every maximal ideal $\fm \subset R$, the $G_{\bA^1_{k_\fm}}$-torsor $\sE_{\bA^1_{k_\fm}}$ is trivial;

\m \label{Gille-2}
letting $\fm \subset R$ range over all the maximal ideals, the following map is surjective:
\[
\tst \q G(R\llp t \rrp)/G(R\llb t \rrb) \surjects 
\prod_\fm G(k_\fm\llp t \rrp)/G(k_\fm\llb t \rrb). 
\]
\eenum
For \ref{Gille-1}, since $\sE_{\bA^1_{k_\fm}}$ is trivial away from $Z_{k_\fm}$, we may glue it arbitrarily with the trivial $G_{k_\fm\llb t\rrb}$-torsor to obtain a $G_{\bP^1_{k_\fm}}$-torsor whose pullback along the infinity section is trivial. By \cite{Gil02}*{Lemme~3.12} (see also \cite{Gil05}), such torsors are trivial over $\bA^1_{k_\fm}$, so $\sE_{\bA^1_{k_\fm}}$ is trivial, that is, \ref{Gille-1} holds. 

For the claim \ref{Gille-2}, 
we first note that the isotropicity assumption implies that $G$ has a proper parabolic subgroup $P \subset G$ (see \Cref{def:tot-iso}). 
By \cite{Gil09}*{Fait 4.3,~Lemme~4.5} (this is where we use the assumptions on $G$), the Whitehead group of the base changes of $G$ is of unramified nature, that is,
\[
\tst 
G(k_\fm\llp t \rrp) = 
G(k_\fm\llp t \rrp)^+G(k_\fm\llb t \rrb) \qxq{for every maximal ideal} \fm \subset R,
\]
where $G(k_\fm\llp t \rrp)^+ \subset G(k_\fm\llp t \rrp)$ is the subgroup generated by $(\sR_u(P))(k_\fm\llp t \rrp)$ and $(\sR_u(P^-))(k_\fm\llp t \rrp)$ with $P^- \subset G$ being a parabolic opposite to $P$ in the sense of \cite{SGA3IIInew}*{Expos\'{e}~XXVI, D\'{e}finition~4.3.3, Corollaire~4.3.5 (i)}. To conclude \ref{Gille-2}, it suffices to show that 
the following pullback maps are surjective: 
\[
\tst (\sR_u(P))(R\llp t \rrp) \surjects \prod_\fm 
(\sR_u(P))(k_\fm\llp t \rrp) \qxq{and} (\sR_u(P^-))(R\llp t \rrp) \surjects \prod_\fm 
(\sR_u(P^-))(k_\fm\llp t \rrp).
\]
For this, we combine the surjectivity of the map $R\llp t\rrp \surjects \prod_\fm 
k\llp t \rrp$ with \cite{SGA3IIInew}*{Expos\'{e}~XXVI, Corollaire~2.5}, according to which both $\sR_u(P)$ and $\sR_u(P^-)$ are isomorphic to affine spaces $\bA^d_R$.
\epf

\bremt
One difference between \Cref{prop:A1-trivial} and some of its versions in the literature 
is that we work directly with the $G_{\bA^1_R}$-torsor $\sE$ instead of first glueing it arbitrarily to a $G_{\bP^1_R}$-torsor and then modifying this extension. Ultimately, this is an expository point, but it highlights that in \ref{Gille-2} there is no need to pursue the analogous surjectivity before taking the quotients. 
\eremt


\section{The quasi-split unramified case of the Grothendieck--Serre conjecture} \label{sec:conclusion}



We turn to the quasi-split case of the Grothendieck--Serre conjecture \ref{conj:GS}: we are ready to settle its semilocal version over unramified regular local rings. By choosing $\cO$ to be either $\bZ$, or $\bQ$, or $\bF_p$ for some prime $p$ and $R$ to be local, this recovers the first assertion in \Cref{thm:main-ann}, see \Cref{eg:geom-regular}.

\bthmt \label{thm:main}
For a Dedekind ring $\cO$, a  semilocal, regular, flat $\cO$-algebra $R$ whose $\cO$-fibers are geometrically regular,\footnote{We recall from \cite{SP}*{Definition~\href{https://stacks.math.columbia.edu/tag/0382}{0382}} that a Noetherian algebra over a field $k$ is \emph{geometrically regular} if its base change to every finite purely inseparable (equivalently, to every finitely generated) field extension of $k$ is regular. 
} 
and a quasi-split reductive $R$-group $G$, 
%
 no nontrivial $G$-torsor trivializes over the total ring of fractions $K \ce \Frac(R)$ of $R$, that is,
\[
\Ker(H^1(R, G) \ra H^1(K, G)) = \{*\}.
\]
\ethmt

\bpf
We pass to connected components to assume that $\Spec(R)$ is connected, so that $R$ is a domain and, in particular, $R \neq 0$. We use \Cref{prop:reduce-to-sc} to replace $G$ by $(G^\der)^{\mathrm{sc}}$, and hence reduce to the case when our quasi-split $G$ is also semisimple and simply connected. 

Let $E$ be a $G$-torsor that trivializes over $K$, so also over $R[\f1r]$ for some $r \in R \setminus \{ 0\}$.  
By Popescu theorem \cite{SP}*{Theorem~\href{https://stacks.math.columbia.edu/tag/07GC}{07GC}}, the ring $R$ is a filtered direct limit of smooth $\cO$-algebras. Thus, a limit argument allows us to assume that $R$ is the semilocalization of a smooth $\cO$-algebra at finitely many primes. 
In this case, \Cref{cor:only-A1} gives a $G_{\bA^1_R}$-torsor $\sE$ whose pullback along the zero section is $E$  such that $\sE$ is trivial away from an $R$-finite closed subscheme $Z \subset \bA^1_R$. By \Cref{prop:A1-trivial} (with \Cref{eg:qsplit}), this $\sE$ is trivial, so $E$ is also trivial, as desired. 
\epf


\begt \label{eg:geom-regular}
In the case when $\cO$ is a perfect field, such as $\bQ$ or $\bF_p$, any regular $\cO$-algebra is geometrically regular, so, for quasi-split $G$, \Cref{thm:main} simultaneously reproves the equicharacteristic case of the Grothendieck--Serre conjecture settled in \cite{FP15} and \cite{Pan20a}. Similarly, in the case when $\cO = \bZ$, the $\bZ$-fibers of a $\bZ$-flat $R$ are geometrically regular if and only if for every prime $p$ and every maximal ideal $\fm \subset R$ of residue characteristic $p$, we have $p \not \in \fm^2$, equivalently, $p$ is a regular parameter for the regular local ring $R_\fm$. In particular, \Cref{thm:main} recovers \Cref{thm:main-ann}.
\eegt

With our main result in hand, we are ready to settle the second assertion of \Cref{thm:main-ann}.

\bthmt \label{thm:split}
For a Dedekind ring $\cO$ and a  semilocal, regular, flat $\cO$-algebra $R$ whose $\cO$-fibers are geometrically regular, 
a reductive $R$-group $G$ is split if and only if so its generic fiber $G_{\Frac(R)}$. 
\ethmt

\bpf
We pass to connected components to assume that $\Spec(R)$ is connected, so that $R$ is a domain, and we set $K \ce \Frac(R)$. Only the `if' part requires an argument, so we assume that $G_{K}$ is split. The geometric fibers of $G$ have a constant type (see \cite{SGA3IIInew}*{Expos\'{e}~XXII, D\'{e}finition~1.13}), and we let $\bbG$ be a split reductive $R$-group of this type, so that $G$ is a form of $\bbG$ that corresponds to some $x \in H^1(R, \underline{\Aut}(\bbG))$ whose pullback to $H^1(K, \underline{\Aut}(\bbG))$ is trivial. We wish to show that $x$ is trivial.

By \cite{SGA3IIInew}*{Expos\'{e}~XXIV, Th\'{e}or\`{e}me~1.3}, we have a short exact sequence of group schemes
\[
1 \ra \bbG^\ad \ra  \underline{\Aut}(\bbG) \ra  \underline{\Out}(\bbG) \ra 1
\]
that, via a fixed pinning of $\bbG$, is split by a homomorphism 
\[
\underline{\Out}(\bbG) \hra \underline{\Aut}(\bbG),
\]
whose source is a constant $R$-group. Any $\underline{\Out}(\bbG)$-torsor $E$ is constant \'{e}tale locally on $R$, so, by \cite{SGA3II}*{Expos\'{e}~X, Corollaire~5.14}, its connected components are finite \'{e}tale over $R$. Thus, by, for instance, \cite{modular-description}*{Lemma~3.1.9}, every $K$-point of $E$ extends to an $R$-point, to the effect that no nontrivial $\underline{\Out}(\bbG)$-torsor trivializes over $K$.

The nonabelian cohomology exact sequence now lifts $x$ to an $\wt{x} \in H^1(R, \bbG^\ad)$ and, since the map 
\[
\underline{\Aut}(\bbG)(K) \surjects \underline{\Out}(\bbG)(K)
\]
is surjective due to the splitting, it also shows that the pullback of $\wt{x}$ to $H^1(K, \bbG^\ad)$ is trivial. \Cref{thm:main} then implies that $\wt{x}$ itself is trivial, and then so is $x$. 
\epf

The ideas of the preceding proof also give a version for quasi-split groups in \Cref{thm:quasi-split}. To put it into context, we recall the following conjecture, which may be traced to results of \cite{CT79} or \cite{Pan09}. Even though not formulated there explicitly, it is sometimes attributed to Colliot-Th\'{e}l\`{e}ne or Panin.

\bconjt
For a regular local ring $R$, if the generic fiber of a reductive $R$-group scheme $G$ has a parabolic subgroup, then $G$ itself has a parabolic subgroup of the same type. 
\econjt

This conjecture ``of Grothendieck--Serre type'' seems to lie deeper than the Grothendieck--Serre conjecture: even in equicharacteristic, it is only known in few cases, see \cite{CT79}, \cite{Pan09}, \cite{PP10}, \cite{PP15}, \cite{Scu18} for precise results. We use the ideas of this article to settle its equicharacteristic case for minimal parabolics, that is, for Borel subgroups, as follows.

\bthmt \label{thm:quasi-split}
Let $R$ a regular semilocal ring, let $K \ce \Frac(R)$ be its total ring of fractions, and let $G$ be a reductive $R$-group such that every form $\cG$ of $G^\ad$ satisfies 
\[
H^1(R, \cG) \hra H^1(K, \cG)
\]
\up{this holds for every $G$ if $R$ contains a field}. Then $G$ is quasi-split if and only if $G_{K}$ is quasi-split. 
\ethmt

\bpf
The injectivity assumption is a special case of the Grothendieck--Serre conjecture and of the ``change of origin'' twisting bijections in nonabelian cohomology \cite{Gir71}*{Chapitre~III, Proposition~2.6.1~(i)}, so the parenthetical assertion follows from the known equicharacteristic case of the Grothendieck--Serre conjecture, see \S\ref{pp:known-cases}. By \cite{Guo20}*{Proposition~14} (whose proof is similar to that of \Cref{thm:split} above), this assumption implies that $G$ is the unique reductive model of its generic fiber, so all we need to do is to assume that $G_K$ is quasi-split and to produce a quasi-split reductive $R$-model of $G_K$. By the properness of the scheme of Borel subgroups, there is an open subscheme $U \subset \Spec(R)$ whose complement is of codimension $\ge 2$ such that even $G_U$ has a Borel subgroup.

Analogously to the proof of \Cref{thm:split}, we reduce to the setting when $\Spec(R)$ is connected, we have a split reductive $R$-group $\bbG$, and $G$ corresponds to an element $x \in H^1(R, \underline{\Aut}(\bbG))$. We fix a Borel subgroup $\bbB \subset \bbG$ that arises from a pinning of $\bbG$, and we consider the subfunctor 
\[
\underline{\Aut}(\bbG, \bbB) \subset \underline{\Aut}(\bbG)
\]
that parametrizes those automorphisms that preserve $\bbB$.  In a reductive group, any two Borels are Zariski locally conjugate, so we are reduced to showing that for our $x \in H^1(R, \underline{\Aut}(\bbG))$ such that $x|_{U}$ lifts to $H^1(U, \underline{\Aut}(\bbG, \bbB))$, the restriction $x|_K \in H^1(K, \underline{\Aut}(\bbG))$ lifts to $H^1(R, \underline{\Aut}(\bbG, \bbB))$.

By \cite{SGA3IIInew}*{Expos\'{e}~XXIV, Th\'{e}or\`{e}me~1.3, Proposition~2.1}, letting $\bbB^\ad \subset \bbG^\ad$ be the Borel subgroup of $\bbG^\ad$ corresponding to $\bbB$, we have a morphism of short exact sequences of group schemes
\[
\xymatrix@R=12pt{
1 \ar[r] & \bbB^\ad \ar@{^(->}[d] \ar[r] & \underline{\Aut}(\bbG, \bbB) \ar@{^(->}[d] \ar[r] & \underline{\Out}(\bbG) \ar@{=}[d] \ar[r] & 1 \\
1 \ar[r] & \bbG^\ad \ar[r] & \underline{\Aut}(\bbG) \ar[r] & \underline{\Out}(\bbG) \ar[r] & 1
}
\]
that, due to our fixed pinning, are compatibly split by some homomorphism 
\[
\underline{\Out}(\bbG) \hra \underline{\Aut}(\bbG, \bbB).
\]
We may first map $x$ to an $\ov{x} \in H^1(R, \underline{\Out}(\bbG))$ and then map $\ov{x}$ via the splitting to obtain a 
\[
y \in H^1(R, \underline{\Out}(\bbG, \bbB))
\]
whose image in $H^1(R, \underline{\Out}(\bbG))$ is also $\ov{x}$. Twisting by (the images of) $y$ gives us the morphism of short exact sequences of $R$-groups of corresponding forms:
\[
\xymatrix@R=12pt{
1 \ar[r] & \cB \ar@{^(->}[d] \ar[r] & \cA_0 \ar@{^(->}[d] \ar[r] & \cE \ar@{=}[d] \ar[r] & 1 \\
1 \ar[r] & \cG \ar[r] & \cA \ar[r] & \cE \ar[r] & 1
}
\]
and, via the ``change of origin'' bijections \cite{Gir71}*{Chapitre~III, Proposition~2.6.1~(i)}, we obtain an $x' \in H^1(R, \cA)$ such that $x'|_{U}$ lifts to $H^1(U, \cA_0)$ for which we need to lift $x'|_K \in H^1(K, \cA)$  to $H^1(R, \cA_0)$ or even to $H^1(R, \cB)$. 

By the nonabelian cohomology sequence, $x'|_U$ even lifts to some $b \in H^1(U, \cB)$. By descent, $\cB \subset \cG$ is the inclusion of a Borel $R$-subgroup, and we let 
\[
\cT \ce \cB/\sR_u(\cB)
\]
be the indicated torus. The image of $b$ is a $t \in H^1(U, \cT)$, which, by purity for torsors under tori \cite{CTS79}*{Corollaire~6.9}, extends uniquely to a 
\[
\wt{t} \in H^1(R, \cT).
\]
Any Levi $R$-subgroup of $\cB$ splits the surjection $\cB \surjects \cT$, and then $\wt{t}$ gives a $\wt{b} \in H^1(R, \cB)$ whose image in $H^1(K, \cB)$, thanks to \cite{SGA3IIInew}*{Expos\'{e}~XXVI, Corollaire~2.3}, is nothing else but $b|_K$. In particular, the image of $\wt{b}$ in $H^1(K, \cA)$ is $x'|_K$, to the effect that $\wt{b}$ is the desired lift.
\epf

We thank Uriya First for pointing out the following further consequence about quadratic forms. 

\bcort \label{cor:q-forms}
For a regular semilocal ring $R$ as in Theorem \uref{thm:main} with $2 \in R^\times$, we have
\[
H^1(R, \SO_{n}) \hra H^1(\Frac(R), \SO_{n}) \qxq{and} H^1(R, \mathrm{O}_{n}) \hra H^1(\Frac(R), \mathrm{O}_{n}) \qxq{for all} n \ge 1;
\]
moreover, no two nonisomorphic quadratic forms over $R$ that are nondegenerate \up{in the sense that their associated symmetric bilinear forms are perfect} become isomorphic over $\Frac(R)$. 
\ecort

\bpf
Every inner form of $\SO_n$ is an $\SO(E)$ for a nondegenerate quadratic space $E$ over $R$ of rank $n$. Thus, by twisting \cite{Gir71}*{Chapitre~III, Proposition~2.6.1~(i)}, the injectivity assertion for $\SO_n$ reduces to showing that
\[
\Ker(H^1(R, \SO(E)) \ra H^1(\Frac(R), \SO(E))) = \{*\}.
\]
By the analysis of the long exact cohomology sequence \cite{CT79}*{page~17, proof of (D)$\Leftrightarrow$(E)}, this triviality of the kernel is, in turn, equivalent to its analogue for $\mathrm{O}(E)$. Thus, by twisting again, we are reduced to the injectivity assertion for $\mathrm{O}_n$, which itself, for varying $n$, is a reformulation of the assertion about quadratic forms. For the latter, however, due to the cancellation theorem for quadratic forms, specifically, due to \cite{CT79}*{Proposition~1.2 (D)$\Leftrightarrow$(F)}, we may assume that one of the forms is a sum of copies of the hyperbolic plane. In terms of $\mathrm{O}_n$-torsors, this means that it suffices to show that
\[
\Ker(H^1(R, \mathrm{O}_n) \ra H^1(\Frac(R), \mathrm{O}_n)) = \{*\} \qxq{for all even} n \ge 1.
\]
We then use \cite{CT79}*{page~17, proof of (D)$\Leftrightarrow$(E)} again to replace $\mathrm{O}_n$ by $\SO_n$ in this display. With this replacement, however, the desired triviality of the kernel is a special case of \Cref{thm:main}. 
\epf




\end{document}